\documentclass[12pt]{amsart}
\usepackage{amsmath,amscd,amssymb,amsfonts}
\def\vers{Dec.~2, 2010, v.3}
\def\msum{\h{$\sum$}}
\def\mopl{\h{$\bigoplus$}}
\def\motim{\h{$\bigotimes$}}
\def\sotim{\h{$\otimes$}}
\def\mprod{\h{$\prod$}}
\def\mwedge{\h{$\bigwedge$}}
\def\nin{\noindent}
\def\bs{\bigskip}
\def\ms{\medskip}
\def\bt{\h{$\boxtimes$}}
\def\bts{\boxtimes}
\def\C{{\mathbf C}}
\def\Db{{\mathbf D}}
\def\Dc{{\mathcal D}}
\def\ep{\varepsilon}
\def\h{\hbox}
\def\H{{\mathcal H}}
\def\I{{\mathcal I}}
\def\Kb{K^{\sb}}
\def\Kc{{\mathcal K}}
\def\Kct{\widetilde{\mathcal K}}
\def\Lc{{\mathcal L}}
\def\la{\lambda}
\def\Mc{{\mathcal M}}
\def\Mb{M^{\sb}}
\def\N{{\mathbf N}}
\def\Nc{{\mathcal N}}
\def\Oc{{\mathcal O}}
\def\Omt{\widetilde{\Omega}}
\def\Omu{\underline{\underline{\Omega}}}
\def\Q{{\mathbf Q}}
\def\q{\quad}
\def\rd{\partial}
\def\R{{\mathbf R}}
\def\s{\sigma}
\def\Sf{\mathfrak S}
\def\Sc{{\mathcal S}}
\def\Sn{{\mathcal S}^n\!}
\def\Sym{{\rm Sym}}
\def\X{{\mathcal X}}
\def\Xo{\overline{X}}
\def\Xt{\widetilde{X}}
\def\Y{{\mathcal Y}}
\def\Z{{\mathbf Z}}
\def\DR{\h{\rm DR}}
\def\Gr{\h{\rm Gr}}
\def\IC{\h{\rm IC}}
\def\Ker{\h{\rm Ker}}
\def\Diff{{\rm Diff}}
\def\End{\h{\rm End}}
\def\Ext{\h{\rm Ext}}
\def\Hom{\h{\rm Hom}}
\def\cHom{{\mathcal H}om}
\def\MHM{\h{\rm MHM}}
\def\MHS{\h{\rm MHS}}
\def\rat{\h{\rm rat}}

\def\into{\hookrightarrow}
\def\sb{\raise.2ex\h{${\scriptscriptstyle\bullet}$}}
\def\sc{\,\raise.2ex\h{${\scriptstyle\circ}$}\,}
\def\simto{\buildrel\sim\over\longrightarrow}

\begin{document}
\title{Symmetric Products of Mixed Hodge Modules}
\author{Laurentiu Maxim}
\address{Department of Mathematics, University of Wisconsin-Madison,
480 Lincoln Drive, Madison WI 53706-1388 USA}
\email{maxim@math.wisc.edu}
\author{Morihiko Saito}
\address{RIMS Kyoto University, Kyoto 606-8502 Japan}
\email{msaito@kurims.kyoto-u.ac.jp}
\author{J\"org Sch\"urmann}
\address{Mathematische Institut, Universit\"at M\"unster,
Einsteinstr. 62, 48149 M\"unster, Germany}
\email{jschuerm@math.uni-muenster.de}
\date{\vers}
\begin{abstract}
Generalizing a theorem of Macdonald, we show a formula for the
mixed Hodge structure on the cohomology of the symmetric products
of bounded complexes of mixed Hodge modules by showing the existence
of the canonical action of the symmetric group on the multiple
external self-products of complexes of mixed Hodge modules.
We also generalize a theorem of Hirzebruch and Zagier on the
signature of the symmetric products of manifolds to the case of
the symmetric products of symmetric parings on bounded complexes
with constructible cohomology sheaves where the pairing is not
assumed to be non-degenerate.\end{abstract}

\maketitle

\centerline{\bf Introduction}

\bs\nin
For a complex algebraic variety $X$, let $\Sn X$ denote the
$n$-fold symmetric product of $X$.
This is by definition the quotient of the $n$-fold self-product
$X^n$ by the action of the symmetric group $\Sf_n$.
We assume $X$ is quasi-projective so that $\Sn X$ is an algebraic
variety. Let $\pi:X^n\to\Sn X$ denote the canonical projection.
We have a canonical isomorphism of sheaves of $\Q$-vector spaces
$$\Q_{\Sn X}=(\pi_*\Q_{X^n})^{\Sf_n},$$
where the right-hand side is the $\Sf_n$-invariant part.
This implies canonical isomorphisms
$$H^{\sb}(\Sn X,\Q)=H^{\sb}(X^n,\Q)^{\Sf_n}=
\bigl(\motim^nH^{\sb}(X,\Q)\bigr){}^{\Sf_n},$$
where the last isomorphism follows from the multiple K\"unneth
formula.
So we get an isomorphism of bigraded vector spaces
$$H^{\sb}(\Sc^{\sb}\!X,\Q)\,(:=\mopl_{n\ge 0}\,H^{\sb}(\Sn X,\Q))
=\Sym^{\sb}\!H^{\sb}_{\rm even}(X,\Q)\otimes\mwedge^{\sb}\!
H^{\sb}_{\rm odd}(X,\Q),$$
where $H^{\sb}_{\rm even}(X,\Q):=\mopl_{j:\,{\rm even}}\,H^j(X,\Q)$
(similarly for $H^{\sb}_{\rm odd}$), and $\Sym^{\sb}V{}^{\sb}$ on
the right-hand side denotes the direct sum of the usual symmetric
tensor products $\mopl_{n\ge 0}\,\Sym^nV^{\sb}$ for a graded
vector space $V^{\sb}$.
The last formula, which was implicit in [MS], is
pointed out to us by S.~Kimura, and is closely related to the theory
of finite dimensional motives [Ki] (see also [dBa]).
We show that the above isomorphisms are compatible with mixed Hodge
structures, and extend to the case of intersection cohomology
(with local system coefficients).
These imply formulas in [MS] for the Hodge numbers and Hirzebruch's
$\chi_y$-genus [Hi] of the (intersection) cohomology of symmetric
products, generalizing earlier work by [Ma], [Ch] and others.

In the case when $X$ is smooth, $\Sn X$ is a complex $V$-manifold,
hence is a $\Q$-homology manifold, so the intersection
cohomology of $\Sn X$ coincides with its usual cohomology (see [GM]).
If furthermore $X$ is projective (or, more generally, compact with
a K\"{a}hler desingularization), then the pure Hodge structure on
the cohomology of $X$ defined in [St] (which is reproduced in [PS],
Section 2.5) coincides with the mixed Hodge
structure in [D3] and also with the pure Hodge structure on the
intersection cohomology in [Sa1] (which is obtained by applying the
decomposition theorem to the desingularization in the
non-projective case), see Proposition~(2.8).
Here we can prove only a weaker version of [DB], Th.~5.3 showing
that if an algebraic complex $V$-manifold $X$ is globally embeddable
into a smooth variety (e.g. if $X$ is quasi-projective), then the
filtered Steenbrink complex $(\Omt^{\sb}_X,F)$ is canonically
isomorphic to the filtered Du Bois complex $(\Omu^{\sb}_X,F)$ as
filtered differential complexes on the {\it ambient} variety.
It is, however, unclear whether the isomorphism holds as filtered
differential complexes on the original variety $X$ as noted in
loc.~cit., see also a remark after (2.7.2) below.

In this paper we extend the above assertions on the symmetric
products to the case of arbitrary bounded complexes of mixed Hodge
modules $\Mc\in D^b\MHM(X)$, where $\MHM(X)$ is the abelian category
of (algebraic) mixed Hodge modules on $X$, and $D^b\MHM(X)$ is the
derived category of bounded complexes of $\MHM(X)$.
There are, however, certain technical difficulties associated to
mixed Hodge modules.

For instance, it is not clear a priori if there is a canonical action
of the symmetric group $\Sf_n$ on the $n$-fold external self-product
$\bt^n\Mc$ in a compatible way with the natural action on the
underlying $\Q$-complexes, since the difference in the $t$-structures
of the underlying $\Dc$-modules and $\Q$-complexes gives certain
differences of signs.
In this paper we solve this problem by showing a cancellation of
sings appearing in the morphisms of a commutative diagram,
see Prop.~(1.5) and Th.~(1.9) below.
It is rather surprising that the sign coming from certain changes
of orders of the multiple external products of complexes cancels out
with the sign coming from the anti-commutativity of the exterior
algebra $\bigwedge^{\sb}\Theta$ where $\Theta$ is the sheaf of
vector fields.

We also prove the multiple K\"unneth formula for the $n$-fold
external products of bounded complexes of mixed Hodge modules in a
compatible way with the action of the symmetric group $\Sf_n$,
see (1.12).
For the compatibility with the action of $\Sf_n$, we have to
construct a canonical multiple K\"unneth isomorphism.
Once a canonical morphism is constructed, the assertion is reduced to
the formula for the underlying $\Q$-complexes, which is well known.

As a consequence of these considerations, we get the following
assertion, which is used in [MS]:

\ms\nin
{\bf Theorem~1.} {\it For any bounded complex of mixed Hodge modules
$\Mc$ on a complex quasi-projective variety $X$, the symmetric
product can be defined by
$$\Sn\Mc=(\pi_*\bt^n\Mc)^{\Sf_n}\in D^b\MHM(\Sn X),$$
and we have canonical isomorphisms of graded mixed Hodge structures
$$H^{\sb}(\Sn X,\Sn\Mc)=H^{\sb}(X^n,\bt^n\Mc)^{\Sf_n}
=\bigl(\motim^nH^{\sb}(X,\Mc)\bigr){}^{\Sf_n},$$
in a compatible way with the corresponding isomorphisms
of the underlying $\Q$-complexes.}

\ms
Forgetting mixed Hodge structures, the last assertion of Theorem~1
also holds for any bounded $A$-complexes with constructible
cohomology sheaves $K\in D_c^b(X,A)$ on a topological stratified
space $X$ as in (3.1) with $\dim_AH^{\sb}(X,K)<\infty$ where
$A$ is a field of characteristic 0, see (3.8) below.

In Theorem~1 we use the splitting of an idempotent in
$D^b\MHM(\Sn X)$ (see [BS] and also [LC] for a simpler argument)
together with the complete reductivity of $\Sf_n$ to show the
existence of the $\Sf_n$-invariant part in $D^b\MHM(\Sn X)$.
The complete reductivity is also used to show the commutativity of
the $\Sf_n$-invariant part with the direct image by $X^n\to pt$.
Note that a certain amount of representation theory is needed to
justify the definition of the projector
$\frac{1}{n!}\sum_{\s\in\Sf_n}\s$ defining the invariant part in case
of bounded complexes of mixed Hodge modules,
see Remark~(2.6) below.

\ms\nin
{\bf Corollary~1.} {\it With the above notation, there is a
canonical isomorphism of bigraded mixed Hodge structures
$$H^{\sb}(\Sc^{\sb}\!X,\Sc^{\sb}\!\Mc)=\Sym^{\sb}H^{\sb}_{\rm even}
(X,\Mc)\otimes\mwedge^{\sb}\!H^{\sb}_{\rm odd}(X,\Mc),$$
where $H^{\sb}(\Sc^{\sb}\!X,\Sc^{\sb}\!\Mc)$ has the bigrading
$\mopl_{j,n\ge 0}\,H^j(\Sc^n\!X,\Sc^n\!\Mc)$, and $\Sym^{\sb}$ on the
right-hand side denotes the direct sum of the usual symmetric tensor
products of the graded mixed Hodge structure
$H^{\sb}_{\rm even}(X,\Mc)$ consisting of the even degree part.}

\ms
For this we use the fact that the $\Sf_n$-invariant part
coincides with the maximal quotient on which the action of $\Sf_n$
is trivial (this follows from the complete reductivity of $\Sf_n$).
Corollary~1 implies Th.~1.1 in [MS] for any bounded complexes of
mixed Hodge module $\Mc$.
More precisely, the Hodge numbers of $\Mc\in D^b\MHM(X)$ are defined
for $p,q,k\in\Z$ by
$$h^{p,q,k}(\Mc):=h^{p,q}(H^k(X,\Mc)):=
\dim_{\C}(\Gr^p_F\Gr^W_{p+q}H^k(X,\Mc)_{\C}),$$
where $H^k(X,\Mc)_{\C}$ denotes the underlying $\C$-vector space
of a mixed Hodge structure.
Taking the alternating sums over $k$, we get the $E$-polynomial
in $\Z[y^{\pm 1},x^{\pm 1}]$:
$$e(\Mc)(y,x):=\sum_{p,q} e^{p,q}(\Mc)\,y^px^q\q
\h{with}\q e^{p,q}(\Mc):=\sum_k (-1)^k h^{p,q,k}(\Mc).$$
For the generating series of the above numbers and polynomials,
Corollary~1 implies the following assertion in [MS], Th.~1.1:

\ms\nin
{\bf Corollary~2.} {\it For any bounded complex of mixed Hodge
modules $\Mc$ on a complex quasi-projective variety $X$, we have the
following identities$\,:$}
$$\aligned &\sum_{n \geq 0} \Bigl(\sum_{p,q,k} h^{p,q,k}(\Sn\Mc)\,
y^px^q(-z)^k\Bigr)\,t^n=\prod_{p,q,k}\Bigl(\frac{1}{1-y^px^qz^kt}
\Bigr)^{(-1)^k h^{p,q,k}(\Mc)},\\
&\sum_{n\geq 0}e(\Sn\Mc)(y,x)\,t^n=\prod_{p,q}\Bigl(
\frac{1}{1-y^px^qt}
\Bigr)^{ e^{p,q}(\Mc)}= \exp\Bigl(\sum_{r \geq 1}e(\Mc)(y^r,x^r)\,
\frac{t^r}{r} \Bigr).\endaligned$$

\ms
Indeed, the isomorphism in Corollary~1 holds after passing to the
bigraded quotients $\Gr_F^{\sb}\Gr_{\sb}^W$, and it implies that
the direct sum decomposition
$V^{\sb}=V^{\sb}_{\rm even}\oplus V^{\sb}_{\rm odd}$ with
$V^{\sb}:=H^{\sb}(X,\Mc)$ gives a multiplicative decomposition
of the left-hand side in a compatible way with the one on the
right-hand side defined by the denominators and the numerators.
So the first equality is reduced to the case
$V^{\sb}=V^{\sb}_{\rm even}$ or $V^{\sb}=V^{\sb}_{\rm odd}$,
and easily follows (using a basis of $\Gr_F^{\sb}\Gr_{\sb}^WV^{\sb}$
if necessary).
Here the degree of $z$ corresponds to the degree of the graded
vector space $V^{\sb}$.
The second equality follows from the first by substituting $z=1$,
and the last equality uses the identity
$-\log\bigl(1-t)=\sum_{i\ge 1}t^i/i$.

Note that Corollary~2 in the constant coefficient case with $X$
singular is stated in [Ch], Prop.~1.1 without any comments about the
compatibility with the mixed Hodge structures as in Proposition~(2.2)
below (although [D3] should be used there).
Note also that the formalism in Appendix A of [MS] does not directly
apply to mixed Hodge modules since the external product $\bt$ for
mixed Hodge modules are not defined by using smooth pull-backs
$pr_i^*$ and tensor product $\motim$
(both have certain shifts of degrees in case of mixed Hodge modules,
which cause the problem of sign of the action of $\Sf_n$).
It seems also rather nontrivial whether the isomorphism in the
multiple K\"unneth formula is independent of the order of
$\{1,\dots,n\}$ if one proves the formula by induction on $n$
reducing to the case $n=2$ without assuming the existence of a
canonical isomorphism for multiple products, even though this
independence is crucial to the proof of the compatibility with the
action of $\Sf_n$.
To avoid this problem, we use multiple tensor product in Section 1.

\ms
We define the $\chi_y$-genus of $\Mc\in D^b\MHM(X)$ in
$\Z[y^{\pm 1}]$ by
$$\chi_{-y}(\Mc):=\sum_{p} f^p(\Mc)\, y^p \q \h{with} \q
f^p(\Mc) :=\sum_{k} (-1)^k \dim_{\C}\Gr^p_F H^k(X,\Mc)_{\C}.$$
Only this $\chi_y$-genus has the corresponding characteristic class
version, see [BSY], [Sc2].
A generating series formula for these characteristic classes of
symmetric products is discussed in [CMSSY].
Since
$$f^p(\Mc)=\msum_{q}\,e^{p,q}(\Mc),\q\chi_{-y}(\Mc)=e(\Mc)(y,1),$$
Corollary~2 implies the following assertion in [MS], Cor.~1.2:

\ms\nin
{\bf Corollary~3.} {\it For any bounded complex of mixed Hodge
modules $\Mc$ on a complex quasi-projective variety $X$, we have the
following equalities$\,:$}
$$\sum_{n \geq 0} \chi_{-y}(\Sn\Mc)\,t^n=\prod_{p} \Bigl(
\frac{1}{1-y^pt}\Bigr)^{ f^{p}(\Mc)}=\exp\Bigl(\sum_{r \geq 1}
\chi_{-y^r}(\Mc)\,\frac{t^r}{r}\Bigr).$$

\ms
Replacing cohomology by cohomology with compact supports,
we get $h_{c}^{p,q,k}(\Mc)$, $e_c(\Mc)(y,x)$, etc. instead of
$h^{p,q,k}(\Mc)$, $e(\Mc)(y,x)$, etc. and the assertions also hold
for those numbers and polynomials as is stated in [MS].
Indeed, we can replace $X$ with a compactification
$\Xo$ and $\Mc$ with the zero extension $j_!\Mc$ where
$j:X\to\Xo$ is the inclusion.

We can apply the above formulas, for instance, to the cases where
$\Mc=(a_X)^*\Q$ with $a_X:X\to pt$ the canonical morphism or
$\Mc=(\IC_X\Lc)[-\dim X]$ with $X$ irreducible and $\Lc$ a polarizable
variation of Hodge structure defined on a smooth Zariski-open subset
of $X$, see also [MS] for more examples and applications.

\ms
Setting $y=1$ in the formula of Corollary~3, we get a formula for
the underlying $\Q$-complex $K$.
This formula also holds in case $X$ is a topological stratified
space as in (3.1) and $K\in D_c^b(X,A)$ with
$\dim_AH^{\sb}(X,K)<\infty$ where $A$ is a field of characteristic 0,
see also [MS], Th.~1.4(a).
For such $X$ and $K$ with $A=\R$, we also have a generalization of
[MS], Th.~1.4(c) on the signature of the induced pairing on the
symmetric products as follows:

Assume there is a pairing
$$\phi:K\otimes_{\R}K\to\Db_X:=a_X^!\R,$$
where $a_X:X\to pt$ is the natural morphism.
This induces a pairing
$$\pi_*\bt^n\phi:\pi_*\bt^nK\otimes\pi_*\bt^nK\to\pi_*\bt^n\Db_X=
\pi_*\Db_{X^n}\buildrel{\rm Tr}\over\longrightarrow\Db_{\Sn X},$$
where the last morphism is the trace morphism
${\rm Tr}:\pi_!\pi^!\Db_{\Sn X}\to\Db_{\Sn X}$ associated with the
adjunction for $\pi_!,\pi^!$.
Restricting this self-pairing to $\Sn K$, we get the induced
pairing
$$\Sc^n\phi:\Sn K\otimes\Sn K\to \Db_{\Sn X}.$$
Assume $\dim_{\R}H_c^{\sb}(X,K)<\infty$, and moreover $\phi$ induces
a graded-symmetric self-pairing $\phi_X$ of the graded vector space
$H_c^{\sb}(X,K)$ (in particular, its restriction to the odd degree
part is anti-symmetric).
The last condition is satisfied if $\phi:K\otimes_{\R}K\to\Db_X$ is
symmetric.
We do not assume, however, that $\phi_X$ on $H_c^{\sb}(X,K)$ is
non-degenerate.
Let $\s(\phi)$ be the signature of $\phi_X$ on $H_c^0(X,K)$.
Let $\rho_i$ be the rank of the induced pairing between $H_c^i(X,K)$
and $H_c^{-i}(X,K)$ $(i\in\Z)$.
Set $\chi(\phi)=\sum_i(-1)^i\rho_i$.
This coincides with the Euler characteristic $\chi_c(X,K)$ if
$\phi_X$ is non-degenerate.
Let $\s(\Sc^n\phi)$ be the signature of the induced pairing
on $H_c^0(\Sn X,\Sn K)$.

\ms\nin
{\bf Theorem~2.} {\it With the above notation and assumption, we
have the identity$\,:$}
$$\sum_{n\ge 0}\s(\Sc^n\phi)\,t^n=\frac
{(1+t)^{\frac{\s(\phi)-\chi(\phi)}{2}}}
{(1-t)^{\frac{\s(\phi)+\chi(\phi)}{2}}}$$

\ms
This generalizes a result of Hirzebruch and Zagier [Za]
which is closely related to the Hirzebruch signature theorem,
see also [MS], Th.~1.4(c).
If $X$ is an even-dimensional complex analytic space and if $X$ is a
$\Q$-homology manifold, then we consider $\IC_X\R=\R_X[\dim X]$
instead of $\R_X$, where the dualizing complex $\Db_X$ is given by
$\R_X[2\dim X]$ (here the Tate twist is omitted).
However, this does not cause any problem of sign since the complex
is shifted by an even degree.

Finally, note that Corollary~3 implies Theorem~2 in the case when
$X$ is projective, $K$ underlies a pure $\R$-Hodge module $\Mc$ of
even weight with strict support, and $\phi$ gives a polarization of
$\Mc$, see (3.7) below.
Note also that Theorem~1 and its corollaries hold also for mixed
$\R$-Hodge modules which is defined in the same way as in [Sa1]
[Sa2] using induction on the dimension of the support with
$\Q$-Hodge structure replaced by $\R$-Hodge structure in the
zero-dimensional case.
Here we assume the local monodromies are quasi-unipotent so that the
$V$-filtrations are indexed by $\Q$.
This is different from [Sa3], 1.11.
Note that the proof of (0.10) (i.e.\ Th.~2.2) in loc.~cit. is
still incomplete even now (the problem is very difficult),
and we have to wait until the detailed version of a paper
quoted there will be published.

\ms
We would like to thank S.~Kimura for useful comments about
[MS], Th.~1.1.
The first named author is partially supported by NSF-1005338.
The second named author is partially supported by Kakenhi 21540037.
The third named author is supported by the SFB 878
``groups, geometry and actions''.

\ms
In Section 1 we show the canonical action of the symmetric group
together with the multiple K\"unneth formula for the multiple
external products of bounded complexes of mixed Hodge modules.
In Section 2 we prove Theorem~1 after showing Proposition~(2.2)
which is related to a certain $\la$-structure.
We also prove Proposition~(2.8) in case of complex $V$-manifolds.
In Section 3 we prove Theorem~2 treating only bounded $\R$-complexes
with constructible cohomology sheaves.

\bs\bs
\centerline{\bf 1. Multiple external products and the symmetric group}

\bs\nin
{\bf 1.1.~Multiple external products of complexes.}
Let $\Kb_i$ be bounded complexes of $A$-modules for $i\in[1,n]$,
where $A$ is a field.
We have the $n$-fold tensor complex $\motim_{i=1}^n\Kb_i$ such that
the $j$-th component is given by
$$\mopl_{|p|=j}\,\bigl(\motim_{i=1}^nK_i^{p_i}\bigr),$$
with $|p|:=\msum_{i=1}^np_i$, and the restriction of the
differential to $\motim_{i=1}^nK_i^{p_i}$ is given by
$$\msum_{i=1}^n\,(-1)^{p_1+\cdots+p_{i-1}}d_i,
\leqno(1.1.1)$$
where $d_i$ denotes also the morphism induced by the differential
$d_i$ of $\Kb_i$.

Let $\Kb_i$ be bounded complexes of sheaves of $A$-modules on
topological spaces $X_i$ for $i\in[1,n]$.
We have the $n$-fold external product $\bt_{i=1}^n\Kb_i$ on
$\prod_{i=1}^nX_i$ such that the $j$-th component is given by
$$\mopl_{|p|=j}\,\bigl(\bt_{i=1}^nK_i^{p_i}\bigr),$$
and the differential is given as in (1.1.1).

Note that the $n$-fold external product of sheaves of $A$-modules
$E_i$ is defined by
$$\bt_{i=1}^n\,E_i=E_1\bts\cdots\bts E_n:=
pr_1^{-1}E_1\otimes_A\cdots\otimes_A pr_n^{-1}E_n,$$
with $pr_i$ the projection to the $i$-th factor.
In the case $X=pt$, the $n$-fold tensor product $\motim_{i=1}^nE_i$
is defined by using the universality for multilinear maps
$E_1\times\cdots\times E_n\to E'$.

We will write ${}_A\bt$ instead of $\bt$ when we have to specify $A$.

\ms\nin
{\bf 1.2.~External products of $\Dc$-modules.}
Let $M_i\in M(\Dc_{X_i})$ for $i\in[1,n]$, where $X_i$ is a
complex manifold and $M(\Dc_{X_i})$ denotes the category of
$\Dc_{X_i}$-modules.
Set $\X=\prod_{i=1}^nX_i$.
We have the $n$-fold external product
$${}_{\Oc}\bt_{i=1}^n\,M_i=M_1\bts\cdots\bts M_n\in M(\Dc_{\X}),$$
which is defined by the scalar extension of
$${}_{\C}\bt_{i=1}^n\,M_i:=
pr_1^{-1}M_1\otimes_{\C}\cdots\otimes_{\C}pr_n^{-1}M_n$$
by the inclusion
$$\Oc'_{\X}:=pr_1^{-1}\Oc_{X_1}\otimes_{\C}\cdots\otimes_{\C}
pr_n^{-1}\Oc_{X_n}\into\Oc_{\X}.
\leqno(1.2.1)$$

For bounded complexes of $\Dc_{X_i}$-modules $\Mb_i\in
C^b(\Dc_{X_i})$, we can then define the $n$-fold external product
$${}_{\Oc}\bt_{i=1}^n\,\Mb_i=\Mb_1\bts\cdots\bts \Mb_n\in
C^b(\Dc_{\X}),$$
where the differential is given as in (1.1.1).
This induces the $n$-fold external product in $D^b(\Dc_{\X})$
for $\Mb_i\in D^b(\Dc_{X_i})$.

Note that the above definition also applies to the case of
$\Oc_X$-modules.

\ms\nin
{\bf 1.3.~Action of the symmetric group $\Sf_n$.}
With the notation of (1.2), assume $X_i=X$ so that $\X=X^n$.
Let $E_i$ be sheaves of $A$-modules on $X$.
For $\s\in\Sf_n$, we have a natural isomorphism (without any sign)
$$\s^{\#}:\bt_{i=1}^n\,E_i\simto\s_*\bigl(\bt_{i=1}^n\,E_{\s(i)}
\bigr).
\leqno(1.3.1)$$
Here $\s_*$, $\s^*$ denote the sheaf-theoretic direct image and
pull-back associated to the action of $\s$ so that
$\s_*=(\s^{-1})^*$, and hence
$\s_*pr_i^*=(\s^{-1})^*pr_i^*=pr_{\s(i)}^*$.
We have moreover
$$\s_*\tau^{\#}\sc\s^{\#}=(\s\tau)^{\#},
\leqno(1.3.2)$$
as morphisms
$$\bt_{i=1}^n\,E_i\to\s_*\tau_*\bigl(\bt_{i=1}^n\,E_{\s\tau(i)}
\bigr)=(\s\tau)_*\bigl(\bt_{i=1}^n\,E_{\s\tau(i)}\bigr).$$
(Set $E'_i=E_{\s(i)}$ so that $E'_{\tau(i)}=E_{\s\tau(i)}$.)
Thus we get a contravariant action.
If one prefers a covariant action, it can be defined by
$\s_{\#}=(\s^{\#})^{-1}.$

Let $\Kb_i$ be bounded complexes of sheaves of $A$-modules on $X$.
We have a canonical isomorphism (see [D2])
$$\s^{\#}:\bt_{i=1}^n\,\Kb_i\simto\s_*\bigl(\bt_{i=1}^n\,
\Kb_{\s(i)}\bigr),
\leqno(1.3.3)$$
which is defined for $m_i\in K^{p_i}_i$ by
$$\aligned &\bt_{i=1}^nm_i\mapsto
(-1)^{\nu(\s,p)}\s_*\bigl(\bt_{i=1}^n\,m_{\s(i)}\bigr),\\
&\h{with}\q\nu(\s,p):=\msum_{i<j,\,\s(j)<\s(i)}\,p_ip_j.\endaligned
\leqno(1.3.4)$$
Here $(-1)^{\nu(\s,p)}$ coincides with the sign of the
permutation of $\{i\in[1,n]\mid p_i:\h{odd}\}$ ignoring the $i$
with $p_i$ even.
More precisely, let $1\le i_1<\dots<i_{n'}\le n$ be the integers
with $p_{i_k}$ odd $(k\in[1,n'])$, and $1\le j_1<\dots<j_{n'}\le n$
be the integers such that $\{\s(i_k)\}=\{j_k\}$.
There is $\s'\in\Sf_{n'}$ such that $j_{\s'(k)}=\s(i_k)$,
and $(-1)^{\nu(\s,p)}$ is the sign of $\s'$.

For $\Mb_i\in C^b(\Dc_{X_i})$, we have similarly a canonical
isomorphism in $C^b(\Dc_{X^n})$
$$\s^{\#}:{}_{\Oc}\bt_{i=1}^n\,\Mb_i\simto\s_*(_{\Oc}\bt_{i=1}^n\,
\Mb_{\s(i)}),
\leqno(1.3.5)$$
inducing an isomorphism in $D^b(\Dc_{X^n})$ for
$\Mb_i\in D^b(\Dc_{X_i})$.
Here $\s_*$ may be viewed as the direct image of $\Dc$-modules
since the action of $\s$ is an isomorphism.

\ms\nin
{\bf 1.4.~Compatibility with the de Rham functor.}
With the notation of (1.2), the external product $\bt$ is compatible
with the de Rham functor $\DR$.
This means that there is a canonical morphism for
$\Mb_i\in C^b(\Dc_{X_i})$:
$${}_{\C}\bt_{i=1}^n\DR_{X_i}(\Mb_i)\to\DR_{\X}({}_{\Oc}\bt_{i=1}^n
\Mb_i),
\leqno(1.4.1)$$
where the left-hand side is defined for $\Kb_i=\DR_{X_i}\Mb_i$ as in
(1.1) with $A=\C$.
Here we use right $\Dc$-modules so that the $j$-th component of
$\DR_X(M)$ for a right $\Dc_X$-module $M$ on a complex manifold $X$
is given by $M\otimes_{\Oc}\mwedge^{-j}\Theta_X$ where $\Theta_X$ is
the sheaf of holomorphic vector fields, and the complex is
locally identified with the Koszul complex associated to the
differential operators $\rd/\rd x_i$ if one chooses a local
coordinate system $(x_1,\dots,x_d)$ so that the
$\mwedge^{-j}\Theta_X$ are locally trivialized, see [Sa1].
Note that (1.4.1) is a quasi-isomorphism in the holonomic case
(i.e. if the $\H^k\Mb_i$ are holonomic $\Dc$-modules).

For $m_i\in M^{p_i}_i$ and $\eta_j\in\mwedge^{-q_j}\Theta_{X_j}$,
the morphism (1.4.1) is then given by
$$\aligned &\bt_{i=1}^n(m_i\sotim\eta_i)\mapsto
(-1)^{\nu(p,q)}\bigl(\bt_{i=1}^nm_i\bigr)\sotim
(\mwedge_{i=1}^n\,pr_i^*\eta_i\bigr),\\
&\,\h{with}\q\q\nu(p,q):=\msum_{i>j}\,p_iq_j,\endaligned
\leqno(1.4.2)$$
where $\mwedge_{i=1}^n\,pr_i^*\eta_i:=pr_1^*\eta_1\wedge\cdots\wedge
pr_n^*\eta_n\in\mwedge^{-\Sigma_i q_i}\Theta_{\X}$, and the sign
comes from the same reason as (1.3.4).

\ms\nin
{\bf Proposition~1.5.} {\it With the notation of $(1.3)$, the action
of the symmetric group $\Sf_n$ is compatible with the de Rham functor
$\DR$, i.e. for $\Mb_i\in C^b(\Dc_X)$, there is a commutative diagram
$$\begin{matrix}{}_{\C}\bt_{i=1}^n\DR_X\bigl(\Mb_i\bigr)&
\buildrel{\s^{\#}}\over\longrightarrow&
\s_*\bigl({}_{\C}\bt_{i=1}^n\DR_X\bigl(\Mb_{\s(i)}\bigr)\bigr)\\
&\raise12pt\h{ }\raise-8pt\h{ }&\downarrow\\
\downarrow&&
\s_*\bigl(\DR_{X^n}\bigl({}_{\Oc}\bt_{i=1}^n\Mb_{\s(i)}\bigr)\bigr)\\
&\raise12pt\h{ }\raise-8pt\h{ }&\downarrow\\
\DR_{X^n}\bigl({}_{\Oc}\bt_{i=1}^n\Mb_i\bigr)&
\buildrel{{\rm DR}(\s^{\#}})\over\longrightarrow&
\DR_{X^n}\bigl(\s_*\bigl({}_{\Oc}\bt_{i=1}^n\Mb_{\s(i)}\bigr)\bigr)
\end{matrix}$$
where the horizontal morphisms are induced by $(1.3.3)$ and
$\DR_{X^n}$ of $(1.3.5)$, and the vertical morphisms are induced
by $(1.4.1)$ together with the commutativity of $\DR_{X^n}$ and
$\s_*$.}

\ms\nin
{\it Proof.}
We may assume $\s=(k,k+1)$ for some $k\in[1,n-1]$, i.e. $\s(i)=i$
for $i\notin\{k,k+1\}$ and $\s\ne id$.
Indeed, $\Sf_n$ is generated by such elements and the action is
compatible with the group law by (1.3.2).
(The last property cannot be used to define the action of $\s\in\Sf_n$
without showing the independence of factorizations of $\s$.)

Let $m_i\in M^{p_i}_i$, $\eta_j\in\mwedge^{-q_j}\Theta_X$.
Consider the image of $\bt_{i=1}^n(m_i\sotim\eta_i)$ in each
term of the diagram.
These are given up to sign by
$$\begin{matrix}\bt_{i=1}^n(m_i\sotim\eta_i)&\to&
\s_*\bigl(\bt_{i=1}^n(m_{\s(i)}\sotim\eta_{\s(i)})\bigr)\\
&&\downarrow\\
\downarrow&&\s_*(\bt_{i=1}^nm_{\s(i)})\otimes
\s_*(\mwedge_{i=1}^n\,pr_i^*\eta_{\s(i)})\\
&&\downarrow\\
(\bt_{i=1}^nm_i)\otimes(\mwedge_{i=1}^n\,pr_i^*\eta_i)&\to&
\s_*(\bt_{i=1}^nm_{\s(i)})\otimes(\mwedge_{i=1}^n\,pr_i^*\eta_i)
\end{matrix}$$
We have to show that the signs associated to the morphisms
of the diagram cancel out.
Since $\s(i)=i$ for $i\notin\{k,k+1\}$, certain signs associated
to the two vertical morphisms cancel out.
Indeed, these are associated to the sum of $p_jq_j$ over
$i<j$ with $(i,j)\ne(k,k+1)$ in (1.4.2).
So the assertion is reduced to the case $n=2$ and $\s\ne id$.
Then the signs coming from the horizontal and vertical morphisms
are given by
$$(-1)^{(p_1+q_1)(p_2+q_2)},\q(-1)^{p_1p_2}\q\h{and}\q
(-1)^{p_2q_1},\q(-1)^{p_1q_2},\q(-1)^{q_1q_2},$$
where the last sign comes from the anti-commutativity in
$\mwedge^{-q_1-q_2}\Theta_{X^2}$, and the other signs come
from (1.3.4) and (1.4.2).
So the assertion follows.

\ms\nin
{\bf Remarks~1.6.} (i)
In the above argument it is also possible to use left $\Dc$-modules
instead of right $\Dc$-modules if we replace the de Rham functor
$\DR_X$ with $\DR_X[-\dim X]$ and cancel the effect of the shift
of complexes by using the twist of the character ${\ep}^{\dim X}$.
Here
$$\ep:\Sf_n\to\{-1,1\}
\leqno(1.6.1)$$
is a character such that $\ep(\s)$ is the sign of a permutation
$\s$.
In this case $\bigwedge^{\sb}\Theta_X$ is replaced by
$\Omega_X^{\sb}=\bigwedge^{\sb}\Omega_X^1$.
Note that the shift of complex $[k]$ in general corresponds
to the twist of the action by the character $\ep^k$.

In the case $M=\Oc_X$, we have to twist the action by $\ep^{\dim X}$
since the de Rham functor $\DR_X$ is shifted by $\dim X$ so that
$\DR_X(\Oc_X)=\C_X[\dim X]$.

\ms
(ii) In Proposition~(1.5), we considered the de Rham complexes in
the category of $\C$-complexes.
However, Proposition~(1.5) holds also if we use the filtered de Rham
functor which is defined for bounded complexes of filtered
$\Dc$-modules $(M^{\sb},F)$ and whose value is in the category
$C^b\!F(\Oc_X,\Diff)$ of bounded filtered differential complexes
in [Sa1].
This filtered differential complex version is needed in [CMSSY].
Here the multiple external products $_{\C}\bt$ in the first row of
the diagram in Proposition~(1.5) is replaced by $_{\Oc}\bt$,
since the multiple external product for bounded filtered differential
complexes is defined by using the scalar extension by (1.2.1).
However, we have to define first $_{\C}\bt$ for bounded filtered
differential complexes before applying the scalar extension by
(1.2.1).
So we need the category $C^b\!F(\Oc'_{\X},\Diff)$
consisting of bounded filtered differential complexes of
$\Oc'_{\X}$-modules, see (1.2.1) for $\Oc'_{\X}$.
This category is defined in the same way as in [Sa1] using
$$\Dc'_{\X}:={}_{\C}\bt_{i=1}^n\Dc_{X_i}=\Oc'_{\X}\langle\rd_1,\dots,
\rd_{d'}\rangle,$$
instead of $\Dc_{\X}=\Oc_{\X}\langle\rd_1,\dots,\rd_{d'}\rangle$
where $\rd_i:=\rd/\rd x_i$ for local coordinates $x_1,\dots,x_{d'}$
of $\X$ with $d'=\dim\X$.

To prove the above variant of Proposition~(1.5), we first prove the
commutativity of the original diagram in Proposition~(1.5) by the
same argument as before where the multiple external product
$_{\C}\bt$ in the first row is defined in $C^b\!F(\Oc'_{\X},\Diff)$.
Then we can take the scalar extension of the two terms in the first
row by the morphism (1.2.1) since the other terms are filtered
differential complexes of $\Oc_{\X}$-modules.

\ms\nin
{\bf 1.7.~Action of $\Sf_n$ on mixed Hodge modules.}
Let $\MHM(X)$ be the category of mixed Hodge modules [Sa2], and
$D^b\MHM(X)$ be the derived category of bounded complexes of
$\MHM(X)$.
Here we assume $X$ is a complex manifold or a smooth complex
algebraic variety.
(Since $X$ is assumed quasi-projective in Theorem~1, we may
replace $X$ with a smooth variety containing it.)
In the algebraic case, we use analytic $\Dc$-modules (assuming
the stratifications are algebraic) to simplify the de Rham
functor.

For $\Mc_i\in\MHM(X)$, there is the $n$-fold external product
$$\bt_{i=1}^n\,\Mc_i:=\Mc_1\bts\cdots\bts\Mc_n\in\MHM(X^n),$$
and there is a natural action of $\Sf_n$ on it by Proposition~(1.5).
This implies a natural action of $\Sf_n$ on
$$\bt_{i=1}^n\,\Mc_i\in C^b\MHM(X^n)\q\h{for}\q
\Mc_i\in C^b\MHM(X),$$
and then on $\bt_{i=1}^n\,\Mc_i\in D^b\MHM(X^n)$ for
$\Mc_i\in D^b\MHM(X)$.
By Proposition~(1.5) this action is compatible with the natural
action on the underlying $\Q$-complexes using the faithfulness of
the functor
$$\otimes_{\Q}\C:D^b_c(X,\Q)\to D^b_c(X,\C).$$
Here the faithfulness follows from the well-known formula
$$\Hom(K,K')=H^0(X,\R\cHom(K,K'))\q\h{for}\q K,K'\in D^b_c(X,A).$$
Note also that the composition of the functor rat associating the
underlying $\Q$-complex and the above functor $\otimes_{\Q}\C$ is
canonically isomorphic to the de Rham functor
$$\DR_X:D^b\MHM(X)\to D^b_c(X,\C).$$
This follows from the construction of the realization functor in
[BBD].

\ms\nin
{\bf 1.8.~Mixed Hodge modules on singular varieties.}
Let $X$ be a complex algebraic variety or a complex analytic space
(assumed Hausdorff).
A filtered $\Dc$-module $(M,F)$ on $X$ is a collection of filtered
$\Dc_Z$-modules $(M_{U\into Z},F)$ for any closed embeddings
$U\into Z$ where $U$ is an open subvariety (or an open subset) of
$X$, $Z$ is smooth, and $(M_{U\into Z},F)\in MF(\Dc_Z)_U$.
Here $MF(\Dc_Z)_U\subset MF(\Dc_Z)$ is defined by the condition
that the $\Gr^F_pM_{U\into Z}$ are $\Oc_U$-modules
(in particular $M_{U\into Z}$ is supported on $U$),
and it is satisfied by mixed Hodge modules supported on $U$,
see [Sa1], Lemma~3.2.6.
They satisfy some compatibility conditions,
see [Sa1], 2.1.20 and also [Sa4], 1.5.
For instance, if there are two closed embeddings
$U_a\into Z_a\,(a=1,2)$, set
$$U_{1,2}:=U_1\cap U_2\into Z_{1,2}:=Z_1\times Z_2.$$
Then we have isomorphisms for $a=1,2$
$$(M_{U_a\into Z_a},F)|_{Z_a\setminus(U_a\setminus U_{a'})}\cong
(pr_{Z_a})_*(M_{U_{1,2}\into Z_{1,2}},F),
\leqno(1.8.1)$$
where $a':=3-a$ and $(pr_{Z_a})_*$ is the the direct image as a
filtered $\Dc$-module under the projection to $Z_a$.
Here we use the following:

Let $f:X'\to Y'$ be a morphism of smooth varieties or complex
manifolds inducing an isomorphism $X\simto Y$ where $X\subset X'$,
$Y\subset Y'$ are closed subvarieties or closed subspaces.
Then the direct image of $\Dc$-modules induces an equivalence of
categories
$$f_*:MF(\Dc_{X'})_X\simto MF(\Dc_{Y'})_Y.
\leqno(1.8.2)$$
This assertion is local, and is reduced to the case where $X\into X'$
is a minimal embedding at $x\in X$, i.e.\ the Zariski tangent space of
$X$ at $x$ has the same dimension as $X'$.

Mixed Hodge modules on singular varieties can be defined in the same
way as above.
In the above notation, we use an equivalence of categories
$$f_*:\MHM(X')_X\simto\MHM(Y')_Y,
\leqno(1.8.3)$$
where $\MHM(X')_X\subset \MHM(X)$ is the full subcategory consisting
of objects supported on $X$.
Note that (1.8.3) holds also for complexes by replacing
$\MHM(X')_X$ with $C^b\MHM(X')_X$ and $\MHM(Y')_Y$ with
$C^b\MHM(Y')_Y$ since
$$f_*\Mc=H^0f_*\Mc,\,\,\,H^jf_*\Mc=0\,(j\ne 0)\q\h{for}\,\,\,
\Mc\in\MHM(X')_X.\leqno(1.8.4)$$

\ms\nin
{\bf Theorem~1.9.} {\it Let $X$ be complex algebraic variety or a
complex analytic space.
Assume $X$ is globally embeddable into a smooth variety or space.
Let $\Mc_i\in C^b\MHM(X)$.
For $\s\in\Sf_n$ we have a contravariant action
$$\s^{\#}:\bt_{i=1}^n\,\Mc_i\simto\s_*(\bt_{i=1}^n\,\Mc_{\s(i)})
\q\h{in}\q C^b\MHM(X^n),$$
satisfying $\s_*\tau^{\#}\sc\s^{\#}=(\s\tau)^{\#}$ for any
$\s,\tau\in\Sf_n$.
The induced isomorphism in $D^b\MHM(X^n)$ for $\Mc_i\in D^b\MHM(X)$
is compatible with the canonical action on the underlying
$\Q$-complexes via the functor $\rat$ associating the underlying
$\Q$-complex.}

\ms\nin
{\it Proof.}
Note first that the assertion is proved in (1.7) if $X$ is smooth.
By hypothesis there is a closed embedding $X\into X'$ with $X'$
smooth.
Then we have an equivalence of categories
$$\MHM(X)\simto\MHM(X')_X,$$
where $\MHM(X')_X\subset \MHM(X')$ denotes the full subcategory
consisting of objects supported on $X$.
This follows from the definition of $\Dc$-modules on complex
varieties or complex analytic spaces, see (1.8).
So any $\Mc\in C^b\MHM(X)$ is canonically represented by
$$\Mc'\in C^b\MHM(X')_X:=C^b\bigl(\MHM(X')_X\bigr).$$
Here $\Mc'$ is a bounded complex of mixed Hodge module defined by
using filtered $\Dc_{X'}$-modules in the usual sense.
Then Theorem~(1.9) in the smooth case can be applied to
$\Mc'$, and we get the canonical action of $\Sf_n$ in
$C^b\MHM((X')^n)_{X^n}$.
This action is independent of the choice of $X'$.
Indeed, if there are two closed embeddings $X\into X'_a\,(a=1,2)$,
set $X'_3=X'_1\times X'_2$.
By (1.8.3) we have an equivalence of categories
$$(pr_a)_*:C^b\MHM(X'_3)_X\simto C^b\MHM(X'_a)_X,$$
induced by the direct image under the projection
$pr_a:X'_3\to X'_a\,(a=1,2)$, and similarly for the projection
between their multiple fiber products where $\MHM(X'_a)_X$ is
replaced by $\MHM((X'_a)^n)_{X^n}$ for $a=1,2,3$.
Since (1.8.4) holds for $f=pr_a$, we have the commutativity of the
multiple external product with the direct image as in (1.12.2) below
also in the analytic case.
So the independence of the embedding follows.
This finishes the proof of Theorem~(1.9).

\ms\nin
{\bf Remark~1.10.} In the case $\Mc_i=\Mc\,(\forall\,i)$ and
$\Mc=a_X^*\Q$ for a variety $X$ or $\Mc=\IC_X\Lc$ with $\Lc$ a
variation of Hodge structure on a smooth dense Zariski-open subset
of an irreducible variety $X$, it is not difficult to show
Theorem~(1.9) using the following property:

\ms\nin
$(P)$ There is a sufficiently small smooth open subset $U$ of $X^n$
which is stable by the action of $\Sf_n$ and such that the
restriction induces an isomorphism
$$\End(\bt^n\!\Mc)\simto\End((\bt^n\!\Mc)|_U).$$
Note that
$$\bt^n a_X^*\Q=a_{X^n}^*\Q,\q\bt^n(\IC_X\Lc)=\IC_{X^n}(\bt^n\Lc).$$

\ms\nin
{\bf 1.11.~Direct image of mixed Hodge modules.}
The definition of the direct image of mixed Hodge modules under
a morphism of complex algebraic varieties $f:X\to Y$ is as follows
(see the proof of [Sa2], Th.~4.3):
Take an affine open covering $U_j$ of $X$, and let
$U_J=\bigcap_{j\in J}U_j$ with $f_J$ the restriction of $f$ to
$U_J$.
Here $X$ may be singular, and we have closed embeddings
$U_j\into Z_j$ with $Z_j$ smooth since $U_j$ are affine,
see (1.8) for mixed Hodge modules on singular varieties.
We may also assume that there is an affine open covering $\{U'_j\}$
of $Y$ together with closed embeddings $U'_j\into Z'_j$ such that
$Z'_j$ is smooth, $f(U_j)\subset U'_j$ and $f|U_j$ is extended
to $f'_j:Z_j\to Z'_j$, see also [Sa1], 2.3.9.
Then, for any bounded complex of mixed Hodge modules $\Mc$ on $X$,
there is a quasi-isomorphisms $\Nc\to\Mc$ such that
$$H^k(f_J)_*(\Nc^p|_{U_J})=0\,\,\,\,\h{for any $p\in\Z$, $k\ne 0$
and $J$.}$$
(This is essentially the same argument as in [Be].)
Here $(f_J)_*$ can be defined by using $f'_J:=\prod_{j\in J}f'_j$.
Combining this construction with the Cech complex for $U_J$,
we get a double complex such that the associated single complex
$\Mc'$ is a representative of $\Mc$ satisfying the condition:
$$H^kf_*(\Mc'{}^p)=0\,\,\,\,\h{for any $p\in\Z$ and $k\ne 0$}.
\leqno(1.11.1)$$
Then the direct image $f_*\Mc$ is defined by $H^0f_*(\Mc'{}^{\sb})$.
This is independent of the choice of the above $\Nc$ by the
standard argument in the theory of derived categories since there
is a quasi-isomorphism $\Nc\to\Mc$ for any $\Mc$.
It is also independent of the choice of $U_j$ by using a refinement
of two affine coverings.

In case $X$ is projective, we can take $U_i$ to be the complement
of a hyperplane section, and the resolution $\Nc\to\Mc$ can be
constructed by using the dual of the Cech complex (using the
direct images with proper supports) which is associated to an
open covering defined by the complements of sufficiently general
hyperplane sections, see [Be].

We can similarly define the direct image with compact support
$f_!$ by the dual argument where the Cech complex associated to the
$U_i$ is replaced by the dual of the Cech complex using the direct
images with proper supports and the directions of the morphisms are
all reversed, e.g. we have $\Mc\to\Nc$ instead of $\Nc\to\Mc$.

\ms\nin
{\bf 1.12.~Multiple K\"unneth formula for mixed Hodge modules.}
For morphisms of complex algebraic varieties $f_i:X_i\to Y_i$,
set $\X=\prod_{i=1}^nX_i$, $\Y=\prod_{i=1}^nY_i$, and
$f=\prod_{i=1}^nf_i:\X\to\Y$.
Let $\Mc_i\in D^b\MHM(X_i)$.
We first show that the direct image commutes with the multiple
external products, i.e. there is a canonical isomorphism
$$\bt_{i=1}^n\,(f_i)_*\Mc_i=f_*(\bt_{i=1}^n\,\Mc_i)\q\h{in}\,\,\,
D^b\MHM(\Y),
\leqno(1.12.1)$$
and this also holds with direct images $(f_i)_*$, $f_*$ replaced
by direct images with proper supports $(f_i)_!$, $f_!$.
Moreover the isomorphism is compatible with the action of $\Sf_n$
in case $X_i=X$ and $Y_i=Y$ for any $i$.
These assertions follow from the definition of the direct image of
bounded complexes of mixed Hodge modules as is explained in (1.11).
We note a proof for the usual direct images.
The argument is similar for the direct images with proper supports.

Take a representative $\Mc'_i$ in (1.11) for each $\Mc_i$ so
that (1.11.1) is satisfied.
Then the isomorphism (1.12.1) follows, since there are natural
isomorphisms of mixed Hodge modules
$$\bt_{i=1}^n\,H^0(f_i)_*\Mc^{\prime p_i}_i=
H^0f_*\bigl(\bt_{i=1}^n\,\Mc^{\prime p_i}_i\bigr).
\leqno(1.12.2)$$
By this argument, (1.12.1) is compatible with the action of
$\Sf_n$ in case $X_i=X$ for any $i$.
Moreover, in case $Y_i=pt$ for any $i$, (1.12.1) is compatible
with the corresponding isomorphism (3.8.1) below for the
underlying $\R$-complexes $K_i$ of $\Mc_i$.
This follows from the definition of the realization functor in
[BBD] using Prop.~3.1.8 in loc.~cit.

In case $Y_i=pt$ for any $i$, we also show that the above proof of
(1.12.1) implies the multiple K\"unneth isomorphism of graded mixed
Hodge structures
$$\motim_{i=1}^n\,H^{\sb}(X_i,\Mc_i)=H^{\sb}(\X,\bt_{i=1}^n\,\Mc_i),
\leqno(1.12.3)$$
and this also holds for cohomology with compact supports.
Moreover, the isomorphism is compatible with the action of $\Sf_n$
in case $X_i$ and $\Mc_i$ are independent of $i$.
Here the category of graded-polarizable mixed $\Q$-Hodge structures
$\MHS$ is naturally identified with $\MHM(pt)$ as in [Sa2]
(where `graded-polarizable' means that the graded quotients of the
weight filtration are polarizable, see [D1]).
The obtained mixed Hodge structure in the constant coefficient case
in [Sa2] coincides with that in [D3], see [Sa5].
(Note that $\Ext^i$ in MHS vanishes for $i>1$ by [Ca] although this
is not used in our argument).

By the above argument, the proof of (1.12.3) is reduced to the
multiple K\"unneth formula for $n$-fold tensor products of bounded
complexes of mixed Hodge structures
$$\motim_{i=1}^n\,H^{\sb}\Nc_i\simto H^{\sb}\bigl(\motim_{i=1}^n\,
\Nc_i\bigr),
\leqno(1.12.4)$$
where $\Nc^p_i:=H^0(f_i)_*\Mc^{\prime p}_i$.
For the proof of (1.12.4), we have a canonical morphism induced
by the natural inclusions
$$\motim_{i=1}^n\,\Ker(d:\Nc_i^{p_i}\to\Nc_i^{p_i+1})\into
\Ker(d:\Nc^p\to\Nc^{p+1}),$$
where $\Nc=\motim_{i=1}^n\,\Nc_i$ and
$p=\sum_{i=1}^np_i$.
This is compatible with the multiple K\"unneth formula for the
$n$-fold tensor products of the underlying complexes of $\Q$-vector
spaces, and the assertion is well known in the latter case.
So (1.12.4) follows.

\ms\nin
{\bf Remark~1.13.}
It is also possible to prove the multiple K\"unneth formula for
mixed Hodge modules by induction on $n$ (reducing to the case $n=2$).
In this case it is not easy to show that the obtained isomorphism is
independent of the choice of the order of $\{1,\dots,n\}$ although
this is essential for the proof of the compatibility with the action
of $\Sf_n$.
In the case $Y=pt$, however, there is a canonical isomorphism for
the underlying $\Q$-vector spaces (see (3.8) below), and
we can use this canonical isomorphism by showing its compatibility
with the mixed Hodge structure.

\bs\bs
\centerline{\bf 2. Symmetric products}

\bs\nin
{\bf 2.1.~Representations of $\Sf_n$.}
Since $\Sf_n$ is a finite group, it is completely reductive, and
every finite dimensional representation of $\Sf_n$ over $\Q$ is
semisimple.
(In fact, this is easily shown by taking a positive definite
symmetric pairing $\langle u,v \rangle$ for any finite dimensional
representation on a $\Q$-vector space $V$ and replacing it with
$\sum_{\s\in\Sf_n}\langle\s(u),\s(v)\rangle$ so that
$\langle\s(u),\s(v)\rangle=\langle u,v \rangle$ for any
$u,v\in V$ and $\s\in\Sf_n$.)

Applying this to the group ring $\Q[\Sf_n]$ viewed as a left
$\Q[\Sf_n]$-module, we get the decomposition by irreducible
characters
$$\Q[\Sf_n]=\mopl_{\chi}\,\Q[\Sf_n]_{\chi},$$
where $\chi$ runs over the irreducible characters of $\Sf_n$,
and $\Q[\Sf_n]_{\chi}$ is the sum of simple left
$\Q[\Sf_n]$-submodules of $\Q[\Sf_n]$ with character $\chi$.
Since this decomposition is compatible with the right action of
$\Q[\Sf_n]$, the direct factors $\Q[\Sf_n]_{\chi}$ are two-sided
ideals of $\Q[\Sf_n]$. Hence
$$\Q[\Sf_n]_{\chi}\,\Q[\Sf_n]_{\chi'}\subset
\Q[\Sf_n]_{\chi}\cap\Q[\Sf_n]_{\chi'}=0\,\,\,\,\h{if}\,\,\,
\chi\ne\chi'.$$
By the above decomposition, there are unique elements
$$e_{\chi}\in\Q[\Sf_n]_{\chi}\,\,\,\,\h{with}\,\,\,
\msum_{\chi}\,e_{\chi}=1\in\Q[\Sf_n].$$
The above property implies that the $e_{\chi}$ are mutually
orthogonal idempotents and $e_{\chi}$ is the identity of
$\Q[\Sf_n]_{\chi}$.

Let $V_{\chi}$ be an irreducible representation over $\Q$ with
character $\chi$.
In the case of symmetric groups, it is known that any
irreducible representation over $\C$ is defined over $\Q$.
So the multiplicity of $V_{\chi}$ in $\Q[\Sf_n]$ as a left
$\Q[\Sf_n]$-module coincides with $\dim V_{\chi}$ by the well-known
orthonormal relation between the irreducible characters $\chi$.
This implies that $\dim \Q[\Sf_n]_{\chi}=(\dim V_{\chi})^2$, and
hence $\Q[\Sf_n]_{\chi}$ is isomorphic to the full endomorphism
algebra $\End_{\Q}(V_{\chi})$.
The last assertion is used in [D6].

As a corollary of the semisimplicity, the $\Sf_n$-invariant part
of a finite dimensional representation on a $\Q$-vector space
is identified with the $\Sf_n$-coinvariant part, which is by
definition the maximal quotient on which the action of $\Sf_n$
is trivial.

By the theory of Young diagrams, the irreducible characters $\chi$
correspond to the Young diagrams, i.e. the partitions
$\la=\{\la_1,\la_2,\dots\}$ of $n$ such that
$\la_1\ge\la_2\ge\cdots$ and $\sum_i\la_i=n$.
The corresponding representation can be constructed by using the
Young symmetrizer.
The trivial character 1 corresponds to the trivial partition
$\{n\}$, and the corresponding idempotent is given by
$$e_1=\h{$\frac{1}{n!}$}\msum_{\s\in\Sf_n}\s\in\Q[\Sf_n]
\leqno(2.1.1).$$
The sign character $\ep$ in (1.6.1) corresponds to
$\{1,1,\dots\}={}^t\{n\}$, and
$$e_{\ep}=\h{$\frac{1}{n!}$}\msum_{\s\in\Sf_n}\ep(\s)\s\in\Q[\Sf_n].
\leqno(2.1.2)$$
It is not difficult to prove (2.1.1--2) by using the above
definition of $e_{\chi}$ via the left action of $\Q[\Sf_n]$ on
itself by considering the condition:
$\tau(\msum_{\s}\,a_{\s}\s)=\msum_{\s}\,a_{\s}\s$ (or
$\ep(\tau)\msum_{\s}\,a_{\s}\s$) for any $\tau\in\Sf_n$,
where $a_{\s}\in\Q$.

In general the relation between $e_{\chi}$ and the Young symmetrizer
is not so trivial.
For the proof of Theorem~1, the explicit form of the
idempotents is not needed.
We will need rather the coincidence of the $\Sf_n$-invariant and
coinvariant part as explained above.

\ms
The following is known to the specialists in a more general
situation (i.e.\ for any decompositions $V^{\sb}=V'{}^{\sb}\oplus
V''{}^{\sb}$), see e.g.\ [D6], Sect.~1, [He], 4.2.
It is closely related to a $\lambda$-ring structure of the
Grothendieck group of a graded vector space
although a canonical isomorphism as graded vector spaces is
finer than an equality in the Grothendieck group.

\ms\nin
{\bf Proposition~2.2.} {\it Let $V^{\sb}$ be a finite dimensional
graded vector space. We have the decomposition
$V^{\sb}=V^{\sb}_{\rm even}\oplus V^{\sb}_{\rm odd}$ by the parity
of the degree.
Then we have a canonical isomorphism of bigraded vector spaces
$$\aligned \mopl_{n\ge 0}\,\bigl(\motim^nV^{\sb}\bigr){}^{\Sf_n}&=
\Sym^{\sb}V^{\sb}_{\rm even}\otimes\mwedge^{\sb}V^{\sb}_{\rm odd},\\
\h{i.e.}\q\bigl(\motim^nV^{\sb}\bigr){}^{\Sf_n}&=\mopl_{n'+n''=n}\,
\bigl(\Sym^{n'}V^{\sb}_{\rm even}\otimes\mwedge^{n''}
V^{\sb}_{\rm odd}\bigr),\endaligned
\leqno(2.2.1)$$
where the action of $\Sf_n$ is defined by identifying the graded
vector space $V^{\sb}$ with a complex with zero differential so
that the sign appears as in $(1.3.4)$ with $X=pt$.
If $V^{\sb}$ is graded mixed Hodge structure, then $(2.2.1)$ is an
isomorphism of bigraded mixed Hodge structures.}

\ms\nin
{\it Proof.} In general, the usual symmetric tensor product
$\Sym^nV'$ of a finite dimensional vector space $V'$ can be
identified with a maximal quotient of $\motim^nV'$ on which the
action of $\Sf_n$ is trivial, see (2.1).
We have a similar assertion for $\mwedge^nV'$ using the action of
$\Sf_n$ twisted by the character $\ep$ in (1.6.1).
Using these, we get a canonical surjection
$$\motim^nV^{\sb}\to\mopl_{n'+n''=n}\,\bigl(\Sym^{n'}
V^{\sb}_{\rm even}\otimes\mwedge^{n''}V^{\sb}_{\rm odd}\bigr)
\leqno(2.2.2)$$
Indeed, for $v:=v_1\otimes\cdots\otimes v_n\in\motim^nV^{\sb}$
with $v_i\in V^{k_i}$, its image is defined by the image of
$$(v_{p_1}\otimes\cdots\otimes v_{p_{n'}})\otimes
(v_{q_1}\otimes\cdots\otimes v_{q_{n''}})$$
in the right-hand side, where $p_i$ and $q_j$ are strictly
increasing sequences such that
$$k_{p_i}\,\h{is even},\q k_{q_j}\,\h{is odd},\q\{p_1,\dots,p_{n'}\}
\,\h{$\coprod$}\,\{q_1,\dots,q_{n''}\}=\{1,\dots,n\},$$
where $n'+n''=n$.
Then the morphism respects the action of $\Sf_n$ where
the action is trivial on the target.
We have moreover a canonical morphism from the right-hand side to
the maximal quotient of $\motim^nV^{\sb}$ divided by the subspace
generated by $\s v-v$ for any $\s\in\Sf_n$ and $v\in\motim^nV^{\sb}$.
So the first assertion follows.

The compatibility with the mixed Hodge structures follows from
the property that any morphism of mixed Hodge structures is strictly
compatible with the mixed Hodge structures [D1] since (2.2.2) is
a morphism of mixed Hodge structures.
This finishes the proof of Proposition~(2.2).

\ms\nin
{\bf 2.3.~Proof of Theorem~1.}
Since $X$ is assumed quasi-projective, we can apply Theorem~(1.9),
and get the canonical action of $\Sf_n$ on
$\pi_*\bt^n\Mc\in D^b\MHM(\Sn X)$ which is compatible with the one
on the underlying $\Q$-complexes.
By the splitting of an idempotent in $D^b\MHM(\Sn X)$ (see [BS]
and also [LC] for a simpler argument) which is applied to $e_1$ in
(2.1.1), we get a direct factor $\Sn\Mc$ of $\pi_*\bt^n\Mc$
endowed with two morphisms in $D^b\MHM(\Sn X)$
$$\Sn\Mc\to\pi_*\bt^n\Mc\to\Sn\Mc,$$
whose composition is the identity.
Note that $\Sn\Mc$ is unique up to a canonical isomorphism using
the above two morphisms together with the forgetful functor
associating the underlying $\Q$-complexes.
So the first assertion follows.

Using the decomposition of $\pi_*\bt^n\Mc$ under the irreducible
characters $\chi$ in (2.1) together with the compatibility of the
direct image functor with the composition of $X^n\to\Sn X\to pt$,
we get the canonical isomorphism
$$H^{\sb}(\Sn X,\Sn\Mc)=H^{\sb}(X^n,\bt^n\Mc)^{\Sf_n}.$$
By the multiple K\"unneth formula (1.12) we get the second canonical
isomorphism
$$H^{\sb}(X^n,\bt^n\Mc)=\motim^nH^{\sb}(X,\Mc).$$
These are compatible with the corresponding isomorphisms for the
underlying $\Q$-complexes.
Then the second isomorphism is compatible with the action of $\Sf_n$
using the action on the underlying $\Q$-complexes.
Thus the remaining assertions are proved.
This finishes the proof of Theorem~1.

\ms\nin
{\bf Remarks~2.4.}
(i) We have
$$\Sn\Mc=a_{\Sn X}^*\Q\in D^b\MHM(\Sn X)\,\,\,\h{if}\,\,\,
\Mc=a_X^*\Q\in D^b\MHM(X),
\leqno(2.4.1)$$
where $a_X:X\to pt$, etc. denote the natural morphisms.
This immediately follows from the characterization of $a_X^*\Q$ in
[Sa2], (4.4.2), i.e. it is uniquely characterized by the conditions
that $\rat(\Mc)=\Q_X$ and $\H^0(X,\Mc)=\Q$ as a mixed Hodge structure.

(ii) We have
$$\Sn\Mc=\IC_{\Sn X}\Sn\Lc\in\MHM(\Sn X)\,\,\,\h{if}\,\,\,
\Mc=\IC_X\Lc\in\MHM(X),
\leqno(2.4.2)$$
where $\Lc$ is a polarizable variation of Hodge structure on a smooth
open subvariety $U$ of an irreducible variety $X$,
and $\Sn\Lc$ is defined on the smooth part of $\Sc^nU$.
This follows from the fact that the intersection complexes are
stable by multiple external product, direct factor, and also by the
direct image $\pi_*$ by a finite morphism $\pi$.
Indeed, the intersection complexes are defined by using the
intermediate direct image, and the latter commutes with the
direct image $\pi_*$ by a finite morphism $\pi$
(since $\pi_*$ is an exact functor of mixed Hodge modules).

\ms\nin
{\bf 2.5.~Proof of Corollary~1.}
The assertion follows from Theorem~1 and Proposition~(2.2).

\ms\nin
{\bf Remark~2.6.} Let $\Mc\in D^b\MHM(X)$ with $K$ its underlying
$\Q$-complex. Set
$$\aligned I&:={\rm Im}(\Q[\Sf_n]\to\End(\pi_*\bt^n\Mc)),\\
I'&:={\rm Im}(\Q[\Sf_n]\to\End(\pi_*\bt^nK)),\\
I''&:={\rm Im}(\Q[\Sf_n]\to\End(\Q^n)),\endaligned$$
where the last representation is given by permutation matrices.
In the notation of (2.1) there are sets of irreducible characters
$\Lambda\supset\Lambda'\supset\Lambda''$ such that
$$I=\mopl_{\chi\in\Lambda}\,\Q[\Sf_n]_{\chi},\q
I'=\mopl_{\chi\in\Lambda'}\,\Q[\Sf_n]_{\chi},\q
I''=\mopl_{\chi\in\Lambda''}\,\Q[\Sf_n]_{\chi},$$
since there are surjections $I\to I'\to I''$.
The last morphism can be defined by restricting to the fiber of $\pi$
over a certain good point of the support of $\pi_*\bt^nK$.

In certain cases (e.g.\ in the constant coefficient case), we have
the equality $I=I'=I''$, and the $\Sf_n$-invariant part is clearly
given by the projector $e_1$ in (2.1.1) looking at the action of
$\Sf_n$ on the fiber over a general point of $\Sn X$.
In this case we would not need the representation theory as is
explained in (2.1).
For a general bounded complex of mixed Hodge modules $\Mc$, however,
it is unclear whether the above three coincide, and we need some
argument as in (2.1).

\ms\nin
{\bf 2.7.~Hodge theory on compact complex $V$-manifolds} [St].
Let $X$ be a complex $V$-manifold with $j:X'\into X$ the
inclusion of the smooth part $X'$ of $X$. Following [St], set
$$\Omt_X^p:=j_*\Omega_{X'}^p.
\leqno(2.7.1)$$
By definition, a complex $V$-manifold $X$ is locally a quotient of a
smooth complex manifold $Y$ by an action of a finite group $G$.
Let $\pi:Y\to X$ denote this quotient morphism
(locally defined on $X$).
Since $\pi$ is finite and $X\setminus X'$ has at least
codimension 2, we have locally a canonical isomorphism
$$\Omt_X^p=(\pi_*\Omega_Y^p)^G,
\leqno(2.7.2)$$
Indeed, this holds on $X'$ and we can apply the Hartogs extension
theorem on $Y$ since the pull-back of
$X\setminus X'$ in $Y$ has at least codimension 2.
(Note that $j_*\Omega_{X'}^{\sb}\ne\R j_*\Omega_{X'}^{\sb}$ even in
the algebraic case by taking the global section functor and applying
[Gr], Th.~1$'$, see the proof of [DB], Th.~5.3.
We can prove only a weaker version of loc.~cit.\ by Proposition~(2.8)
below.)

By Steenbrink [St], there is Hodge theory for compact complex
$V$-manifolds with a K\"ahler desingularization.
(This is reproduced in [PS], Section 2.5.)
As for the relation with the theory of Hodge modules, we first note
the following:

The filtered complex $(\Omt_X^{\sb},F)$, with $F^p$ defined by the
truncation $\s_{\ge p}$ in [D1], is a filtered differential complex
in the sense of [DB], and hence also in the sense of [Sa1], i.e. it
belongs to $C^b\!F(\Oc_X,\Diff)$ in loc.~cit.
This assertion follows from (2.7.2).
Note that the derived category of filtered differential complexes in
[DB] is canonically equivalent to the one in [Sa1] if the variety is
{\it smooth}, see [Fi].

There are canonical isomorphisms (see [BBD], [GM]):
$$\Q_X[\dim X]\simto\IC_X\Q,\q
H^{\sb}(X,\Q)\simto{\rm IH}^{\sb}(X,\Q),
\leqno(2.7.3)$$
since a complex $V$-manifold is a $\Q$-homology manifold,
see Remark~(2.10) below.
Moreover, these isomorphisms are lifted to $\MHM(X)$ and $\MHS$
in the algebraic case, see [Sa2].

For a bounded filtered differential complex
$(L^{\sb},F)\in C^b\!F(\Oc_X,\Diff)$ on a complex manifold, we have
the associated bounded complex of filtered right $\Dc$-modules
$\DR^{-1}_X(L^{\sb},F)$, see [Sa1], 2.2.5.
This is naturally extended to the case of singular spaces
and also to the algebraic case.

Let $(M,F)$ be the underlying filtered right $\Dc$-module of the
polarizable Hodge module $\Mc$ corresponding to the intersection
complex $\IC_X\Q$.
This is represented by filtered right $\Dc$-modules $(M_Z,F)$ for
closed embeddings $U\into Z$ where $U$ is an open subset of $X$ and
$Z$ is smooth, see [Sa1], 2.1.20.
Note that $\H^i\Gr^F_p\DR_Z(M_Z,F)$ is an $\Oc_U$-module by [Sa1],
Lemma~3.2, and is independent of $Z$.
So it is globally well-defined on $X$, and is denoted by
$\H^i\Gr^F_p\DR_X(M,F)$.

\ms\nin
{\bf Proposition~2.8.} {\it Let $X$ be a compact complex
$V$-manifold with a K\"ahler desingularization.
Then the pure Hodge structure on $H^{\sb}(X,\Q)$ in {\rm [St]}
coincides with the one on the intersection cohomology
${\rm IH}^{\sb}(X,\Q)$ which is obtained by using the decomposition
theorem in {\rm [Sa1], [Sa3]} for the desingularization.
Moreover, there is a canonical filtered quasi-isomorphism of
complexes of filtered $\Dc$-modules on $X:$
$$\DR_X^{-1}(\Omt_X^{\sb},F)[\dim X]\simto(M,F),
\leqno(2.8.1)$$
where $(M,F)$ is as above, and we have
$$\h{$\H^i\Gr_F^p\DR_X(M,F)=\Omt_X^p\,$ if $\,i=p-\dim X,\,$ and
$\,0\,$ otherwise.}
\leqno(2.8.2)$$
In case $X$ is algebraic, the pure Hodge structure on $H^{\sb}(X,\Q)$
in {\rm [St]} also coincides with the mixed Hodge structure on
$H^{\sb}(X,\Q)$ in {\rm [D3], [Sa2]}.
If furthermore $X$ is a closed subvariety of a smooth complex
algebraic variety $Z$, then there is an isomorphism in the derived
category of filtered differential complexes on $Z$ in the sense of
{\rm [DB]} or {\rm [Sa1]:}
$$(\Omt_X^{\sb},F)=(\Omu_X^{\sb},F),
\leqno(2.8.3)$$
where $(\Omu_X^{\sb},F)$ is the filtered Du Bois complex in
{\rm [DB]}.}

\ms\nin
{\it Proof.} We first show the analytic case.
Using the decomposition theorem ([Sa3], Th.~0.5) for a K\"ahler
desingularization $\rho:\Xt\to X$, we can show that $(M,F)$ is a
direct factor of $\rho_*(\omega_{\Xt},F)$ (where [Sa1] is enough
in case $\rho$ is projective), and this is compatible with the
$\Q$-structure using Deligne's canonical choice of the
decomposition [D5].
This implies a pure Hodge structure on the intersection cohomology.
By the definition of the direct image of filtered $\Dc$-modules in
the analytic case in [Sa2], 2.13 (applied to $a_X:X\to pt$), it is
then enough to show (2.8.1).

Since the assertion is local on $X$, we may assume that $X$ is a
quotient of $Y$ as above, and moreover, $X$ is a closed analytic
subset of a complex manifold $Z$ so that $(M,F)$ is represented by
a filtered $\Dc_Z$-module $(M_Z,F)$.
Thus the assertion is reduced to showing the canonical filtered
quasi-isomorphism of complexes of filtered $\Dc_Z$-modules
$$\DR_Z^{-1}(\Omt_X^{\sb},F)[\dim X]\simto(M_Z,F).
\leqno(2.8.4)$$

Let $\pi':Y\to Z$ be the composition of $\pi:Y\to X$ and the
inclusion $X\into Z$.
Since $\pi'$ is finite, the direct image  as a filtered right
$\Dc$-module $\pi'_*(\omega_Y,F)$ is a filtered $\Dc_Z$-module and
underlies a pure Hodge module corresponding to an intersection
complex with local system coefficients and with strict support $X$.
So we get
$$(M_Z,F)=(\pi'_*(\omega_Y,F))^G.
\leqno(2.8.5)$$
Indeed, the assertion is clear on a Zariski-open subset $X'$ of $X$
over which $\pi$ is \'etale, and $(M_Z,F)$ is uniquely determined
by its restriction to the complement of $X\setminus X'$ (using [Sa1],
Prop.~3.2.2).

The direct image of filtered differential complexes are defined by
the sheaf-theoretic direct image, and the direct image commutes with
the de Rham functor, see [Sa1], Lemma~2.3.6.
Since $Y$ is smooth, we have
$$\DR_Y(\omega_Y,F)=(\Omega_Y^{\sb},F)[\dim Y].$$
So we get a canonical filtered quasi-isomorphism of complexes of
filtered $\Dc_Z$-modules
$$\DR_Z^{-1}(\pi'_*(\Omega_Y^{\sb},F))[\dim Y]\simto
\pi'_*(\omega_Y,F).$$
This is equivalent to an isomorphism in the derived category
since $\pi'_*(\omega_Y,F)$ is a filtered $\Dc_Z$-module and
$\pi'_*\Omega_Y^p=0$ for $p>\dim Y$.
So the isomorphism is compatible with the action of $G$ since it
is clear on the complement of $X\setminus X'$ where $X'$ is as in
the proof of (2.8.5).
By definition we have
$$\DR_Z^{-1}(\pi'_*\Omega_Y^p):=
\pi'_*\Omega_Y^p\otimes_{\Oc_Z}\Dc_Z,$$
and the action of $G$ is induced by that on $\pi'_*\Omega_Y^p$
(i.e. it is the identity on $\Dc_Z$).
So (2.8.4) follows by taking the $G$-invariant part.

In the algebraic case the mixed Hodge structures on
the cohomology of $X$ defined in [D3] and [Sa2] coincide (see [Sa5])
and the canonical isomorphisms in (2.7.3) are lifted to $\MHM(X)$
and $\MHS$.
So it remains to show (2.8.3).
By the same argument as above, we have the algebraic version of
(2.8.1).
In the case $X$ is globally a closed subvariety of a smooth variety
$Z$, we have a canonical isomorphism in the derived category of
filtered $\Dc_Z$-modules
$$\DR_Z^{-1}(\Omt_X^{\sb},F)[\dim X]=(M_Z,F).$$
By [Sa5], Th.~0.2, we have
$$(M_Z,F)=\DR^{-1}_Z(\Omu_X^{\sb},F)[\dim X].$$
We get thus
$$\DR^{-1}_Z(\Omt_X^{\sb},F)=\DR^{-1}_Z(\Omu_X^{\sb},F).$$
So (2.8.3) follows by using [Sa1], Prop.~2.2.10 and applying the
functor $\DR_Z$ which is denoted by $\widetilde{\DR}_Z$ in loc.~cit.
This finishes the proof of Proposition~(2.8).

\ms\nin
{\bf Remarks~2.9.} (i)
For a complex algebraic $V$-manifold $X$, Proposition~(2.8) implies
the canonical isomorphisms of coherent $\Oc_X$-modules
$$\Omt_X^p=\Gr_F^p\Omu_X^{\sb}[p]\,\,\,\,(p\in\Z),
\leqno(2.9.1)$$
since the assertion is local.
It might be possible to prove [DB], Th.~5.3 by extending the
isomorphisms in (2.9.1) if we have the following vanishing of the
negative extensions in the derived category of filtered differential
complexes in loc.~cit.:
$$\Ext^{p-q+1}\bigl((\Omt_X^p,F),(\Omt_X^q,F)\bigr)=0\q
\h{if}\,\,\,q>p+1.$$
There are, however, no truncations $\tau_{\le k}$ for filtered
differential complexes.
The usual definition in [D1] does not work for filtered
differential complexes even in the filtered acyclic case.

\ms
(ii) In case $X$ is singular it is unclear whether for any filtered
differential complex $(K,F)$, there is a {\it filtered injective
resolution}, i.e.\ a quasi-isomorphism $(K,F)\simto(I,F)$ such that
the $\Gr_F^pI^i$ are injective $\Oc_X$-modules.
If it always exists, then the extension group can be calculated
by using an injective resolution, and a version of Th.~5.3 in [DB]
can be proved where the isomorphism is considered in the derived
category of filtered differential complexes in [Sa1].
(For the derived category in [DB], the definition of homotopy in
loc.~cit.\ is not compatible with the calculation using an injective
resolution.)

\ms
(iii) Let $i:X\into Y$ be a closed embedding of algebraic varieties.
In case $X$ is singular (even if $Y$ is smooth), it is unclear
whether the following direct image functor is fully faithful:
$$i_*:D^b\!F(\Oc_X,\Diff)\to D^b\!F(\Oc_Y,\Diff).$$
For $\Oc_Y$-modules $M$, we have the functor $i_{\Oc}^!$ defined by
$$i_{\Oc}^!M:=\cHom_{\Oc_Y}(\Oc_X,M)=\{m\in M\mid\I_Xm=0\}\subset
M,$$
where $\I_X$ is the sheaf of ideals of $X\subset Y$.
This derives $i^!:D^b_{\rm coh}(\Oc_Y)\to D^b_{\rm coh}(\Oc_X)$ by
using injective resolutions.
However, differential operators do not necessarily preserve
$i_{\Oc}^!M\subset M$ in general,
e.g.\ if $X=\{0\}\subset Y={\rm Spec}\,\C[t]$ with
$M=\C[t,t^{-1}]/\C[t]$ and $i_{\Oc}^!M=\C$.

\ms\nin
{\bf Remark~2.10.}
We say that a complex analytic space $X$ is a $\Q$-homology manifold
if the local cohomology $H^i_{\{x\}}\Q_X$ for any $x\in X$ is
isomorphic to $\Q$ if $i=2\dim X$, and vanishes otherwise.
If a finite group $G$ acts on a complex analytic space $X$ and if
$X$ is a $\Q$-homology manifold, then the quotient complex analytic
space $X/G$ is also a $\Q$-homology manifold,
since the action of the stabilizer $G_x$ of $x$ on the local
cohomology $H^{2\dim X}_{\{x\}}\Q_X$ is trivial.
In particular, a complex $V$-manifold is a $\Q$-homology manifold.

\bs\bs
\centerline{\bf 3. Signature on the symmetric products}

\bs\nin
{\bf 3.1.~Induced pairings.}
Let $X$ be a topological stratified space, i.e. $X$ is a Hausdorff
topological space with a stratification given by an increasing
sequence of closed subspaces $X_k\,(k\ge -1)$ with
$X_{-1}=\emptyset$, $X_d=X\,\,(d\gg 0)$, and for any
$x\in X_d\setminus X_{d-1}$ with $d\ge 0$, there is an open
neighborhood $U_x$ of $x$ in $X$ together with a compact topological
space $L_x$ having an increasing sequence of closed subspaces
$(L_x)_k\,(k\ge -1)$ with $(L_x)_{-1}=\emptyset$ and such that
there is a homeomorphism $U_x\cong\R^d\times C(L_x)$ inducing
$X_k\cap U_x\cong\R^d\times C((L_x)_{k-d-1})$ for any $k\ge d$
(see also [GM] and [Sc1], Def.~4.2.1).
Here $C(Z)$ for a topological space $Z$ denotes the open cone of
$Z$ if $Z\ne\emptyset$, and $C(Z)=pt$ if $Z=\emptyset$.
(We do not assume $X$ equidimensional.)
It is known that the multiple K\"unneth formula holds for bounded
complexes on such spaces having constructible cohomology sheaves and
finite dimensional global cohomology groups (with compact supports),
see (3.8) below.

Let $K\in D_c^b(X,\R)$ endowed with a pairing
$$\phi:K\otimes_{\R}K\to\Db_X,$$
where $\Db_X=a_X^!\R$ (see [Ve]) with $a_X:X\to pt$ the canonical
morphism.
We say that $\phi$ is {\it symmetric} if the composition of the
involution of $K\otimes_{\R}K$ with $\phi$ coincides with $\phi$.

We have the induced pairing $\pi_*\bt^n\phi$ which is the
composition of
$$\aligned \pi_*\bt^nK\otimes\pi_*\bt^nK&\to\pi_*\bigl(\bt^nK
\otimes\bt^nK\bigr)\buildrel{\gamma}\over\cong\pi_*\bt^n
(K\otimes K)\\
&\buildrel{\phi}\over\to\pi_*\bt^n\Db_X=\pi_*\Db_{X^n}
\buildrel{\rm Tr}\over\to\Db_{\Sn X},\endaligned
\leqno(3.1.1)$$
where the last morphism is given by the trace morphism Tr
associated with the adjunction for $\pi_!,\pi^!$.
Note that we have a certain sign for the isomorphism $\gamma$ as
in (1.3.4).

Restricting this self-pairing to $\Sn K$, we get the induced
pairing
$$\Sc^n\phi:\Sn K\otimes\Sn K\to \Db_{\Sn X}.$$
Note that the subcomplex $\Sn K\into\pi_*\bt^nK$ is
defined by using the symmetrizer $e_1$ in (2.1.1).

The above construction is compatible with the global section
functor with compact supports.
Here we assume the nonzero $H_c^i(X,K)$ are bounded and finite
dimensional.
We have the induced self-pairing
$$\phi_X:V^{\sb}\otimes V^{\sb}\to\R\q\h{with}\q V^{\sb}:=
H_c^{\sb}(X,K).$$
This is graded-symmetric if $\phi$ is symmetric.
It induces further
$$\motim^n\phi_X:\motim^n\bigl(V^{\sb}\otimes V^{\sb}\bigr)\to\R.$$
Then the induced pairing $\Sc^n\phi_X$ on $H_c^{\sb}(\Sn X,\Sn K)$
coincides with the restriction of the composition
$$\phi_X^n:=\motim^n\phi_X\sc\gamma':\bigl(\motim^nV^{\sb}\bigr)
\otimes\bigl(\motim^nV^{\sb}\bigr)\buildrel{\gamma'}\over\cong
\motim^n\bigl(V^{\sb}\otimes V^{\sb}\bigr)\to\R,
\leqno(3.1.2)$$
to the $\Sf_n$-invariant part
$$H_c^{\sb}(\Sn X,\Sn K)=\bigl(\motim^nV^{\sb}\bigr){}^{\Sf_n}
\into\motim^nV^{\sb}.
\leqno(3.1.3)$$
Here the last inclusion is defined by using the symmetrizer $e_1$ in
(2.1.1).
Note that we have a certain sign for $\gamma'$ as in (1.3.4).

\ms\nin
{\bf 3.2.~Good bases of the cohomology groups.}
With the notation and the assumption of (3.1), set $r_i:=\dim V^i$.
We have bases $v_{i,1},\dots, v_{i,r_i}$ of $V^i$ for $i\in\Z$
satisfying
$$\phi_X(v_{i,j},v_{i',j'})\ne 0\iff i+i'=0,\,\,j=j'\le\rho_i,
\leqno(3.2.1)$$
where the $\rho_i$ is the rank of $\phi_X$ on
$H_c^i(X,K)\otimes H_c^{-i}(X,K)$.
Note that $\rho_i=\rho_{-i}$ for any i. Set
$$J=\{(i,j)\in\Z\times\N\mid j\le r_i\,(\forall\,i)\}.$$
Let $\Lambda_n\subset\N^J$ consisting of
$\mu=(\mu_{i,j})_{(i,j)\in J}\in\N^J$ satisfying the condition:
$$\msum_{i,j}\,\mu_{i,j}=n,\q\msum_{i,j}\,i\,\mu_{i,j}=0,\q
\mu_{i,j}\in\{0,1\}\,\,\h{for $i$ odd}.
\leqno(3.2.2)$$
Then we have a basis of $(\Sc^nV^{\sb})^0$ defined by the images of
$$v^{\mu}:=\motim_i\bigl(\motim_j\,v_{i,j}^{\mu_{i,j}}\bigr)\bigr)
\,\,\,(\mu\in\Lambda_n),$$
where $\motim_i$ and $\motim_j$ are the ordered tensor products
as in (1.1) which are applied successively. Here we identify
$$\Sc^nV^{\sb}:=\bigl(\motim^nV^{\sb}\bigr){}^{\Sf_n}$$
with the maximal quotient of $\motim^nV^{\sb}$ on which the action
of $\Sf_n$ is trivial.
We have to apply the symmetrizer $e_1$ in (2.1.1) to get an element
in the $\Sf_n$-invariant subspace.

By the above argument, the $v_{i,j}$ for $j>\rho_i$ do not
contribute to the signature of $\Sc^n\phi_X$.
Then, replacing $V^{\sb}$ with the subspace generated by $v_{i,j}$
with $j\le\rho_i$, the proof of Theorem~2 is reduced to the case
where the self-pairing $\phi_X$ on $V^{\sb}$ is non-degenerate.

\ms\nin
{\bf 3.3.~Proof of Theorem~2.}
By the above argument we may assume
$$\rho_i=r_i\,\,(\forall\,i\in\Z).
\leqno(3.3.1)$$
With the notation of of (3.2), let $\iota$ be an involution of
$\Lambda_n$ defined by
$$\iota(\mu)=\mu'\q\h{with}\q\mu'_{i,j}:=\mu_{-i,j}\,\,
(\forall\, i,j).$$
Then
$$\phi_X^n(v^{\mu},v^{\mu'})\ne 0\iff \mu'=\iota(\mu).$$
This gives an orthogonal decomposition of $\Sc^nV^{\sb}$ into
the direct factors of the form
$$V^{\mu}:=\R v^{\mu}+\R v^{\iota(\mu)},$$
which has dimension 1 or 2 depending on whether $\iota(\mu)=\mu$ or
not.
If $\dim V^{\mu}=2$, then this orthogonal direct factor is
hyperbolic, and hence the signature is zero.
So this can be neglected.
Thus it is enough to consider only the $V^{\mu}$ with
$\iota(\mu)=\mu$ (and hence $\dim V^{\mu}=1$).
Set
$$\Lambda'_n=\bigl\{\mu\mid\iota(\mu)=\mu\bigr\}\subset
\Lambda_n.$$
We have an additive structure on $\Lambda':=\coprod_n\Lambda'_n$
defined by
$$(\mu+\nu)_{i,j}=\begin{cases}\mu_{i,j}+\nu_{i,j}&
\h{if $\,\mu_{i,j}+\nu_{i,j}\le 1\,$ for any $(i,j)$ with $i$ odd,}\\
0 &\h{if $\,\mu_{i,j}+\nu_{i,j}>1\,$ for some $(i,j)$ with $i$ odd.}
\end{cases}$$
This additive structure is compatible with an orthogonal direct
sum decomposition
$$V^{\sb}=V^{\sb}_1\oplus V^{\sb}_2,$$
if the latter is compatible with the basis $v_{i,j}$
(i.e. if it it corresponds to a partition of the basis $v_{i,j}$).

The right-hand side of the formula in Theorem~2 is compatible with
the above direct sum decomposition since $\s$ and $\chi$ are additive.
So we first calculate the right-hand side of the formula in the
primitive cases of 2 or 1-dimensional vector subspaces of the form:
$$V'{}^{\sb}=\R\,v_{i,j}+\R\,v_{-i,j}\,\,(i\ne 0)\q\h{or}\q
V'{}^{\sb}=\R\,v_{0,j}.
\leqno(3.3.2)$$
Note that $\sum_{\mu\in\Lambda'}\R v^{\mu}\subset\Sc^nV^{\sb}$ is
generated by the images of the multiple tensor products of
vector subspaces of the form
$\bigl(\motim^{\dim V'}V'{}^{\sb}\bigr){}^0$ where
$V'{}^{\sb}$ is as in (3.3.2).

In the first case of (3.3.2) we have $\sigma_{\phi}=0$ and
$\chi_{\phi}=\pm 2$, depending on the parity of the degree $i$.
So the right-hand side of the formula is given in this case by
$$(1-t^2)^{-1}=1+t^2+t^4+\cdots\q\h{or}\q 1-t^2,$$
depending on the parity of the degree $i$.

In the second case of (3.3.2) we have $\chi_{\phi}=1$ and
$\sigma_{\phi}=\pm 1$, depending on the signature of $\phi_X$.
So the right-hand side of the formula is given in this case by
$$(1-t)^{-1}=1+t+t^2+\cdots\q\h{or}\q(1+t)^{-1}=1-t+t^2-t^3\pm
\cdots,$$
depending on the signature of $\phi_X$.

The compatibility with the above direct sum decomposition is rather
nontrivial for the left-hand side for the odd degree part since
there is a problem of sign associated to $\gamma'$ in (3.1.2).
This is trivial for the even degree elements since they
commute with any elements (even with any odd degree elements)
without any signs.
In the above primitive case with even degrees, we can calculate
the left-hand side of the formula, and verify the formula in these
cases.
So the assertion is proved in the case $V^{\sb}$ has only the even
degree part using the above compatibility with direct sum
decompositions.
Then, by Proposition~(2.2) together with the commutativity of even
degree elements with any elements (without any signs), it now remains
to calculate the left-hand side of the formula in the case $V^{\sb}$
has only the odd degree part.

\ms
So the proof of Theorem~2 is reduced to the calculation in the
next subsection.

\ms\nin
{\bf 3.4.~The odd degree case.}
With the notation and the assumption of (3.3), assume further
$V^{\sb}=V^{\sb}_{\rm odd}$.
Take any $\mu\in\Lambda'_n$ where $\mu_{i,j}=0$ for $i$ even by the
above hypothesis.
We have to calculate the sign of
$$\phi_X^n(e_1(v^{\mu}),e_1(v^{\mu})),$$
see (2.1.1) for $e_1$.
Here we replace $v^{\mu}$ with
$$u:=u_1\sotim u'_1\sotim\,\cdots\,\sotim u_r\sotim u'_r,$$
where $u_k=v_{i_k,j_k}$, $u'_k=v_{-i_k,j_k}$, and $n=2r$.
Then $e_1(u)$ coincides with $e_1(v^{\mu})$ up to the sign
$\ep(\tau)$ of $\tau\in\Sf_n$ such that
$\tau(u)=\ep(\tau)v^{\mu}$.
Note that the action of $\Sf_n$ on $\motim^nV^{\sb}_{\rm odd}$
is twisted by the sign character $\ep$ in (1.6.1).
Since the two signs cancel out, we get
$$\phi_X^n(e_1(v^{\mu}),e_1(v^{\mu}))=
\phi_X^n(e_1(u),e_1(u)).$$

We then replace the second $u$ in
$\phi_X^n(e_1(u),e_1(u))$ by
$$u':=u'_1\sotim u_1\sotim\,\cdots\,\sotim u'_r\sotim u_r.$$
Here we get the first sign $(-1)^r$.
This is the sign of $\tau'$ such that $\tau'(u)=(-1)^ru'$.
Hence
$$e_1(u)=(-1)^re_1(u').$$
By (2.1.1), we have to calculate the sign of
$$\msum_{\s,\s'\in\Sf_n}\,\phi_X^n(\s u,\s'u')=
\msum_{\s\in\Sf_n}\,\phi_X^n(\s u,\s u')=
\phi_X^n(u,u')\,n!.
\leqno(3.4.1)$$
Here the middle equality follows from the vanishing of
$\phi_X^n(\s u,\s'u')$ for $\s\ne\s'$, since $u'$
coincides with $v^{\iota(\mu)}$ up to a sign.
For the last equality, we need
$$\phi_X^n(\s v,\s v')=\phi_X^n(v,v')\q\h{for any}\,\,v,v'\in
\motim^nV^{\sb}_{\rm odd}.
\leqno(3.4.2)$$
This follows from the definition of $\phi_X^n$ in (3.1.2).
Indeed, we have $n=2r$, and the sign of $\gamma'$ in (3.1.2) is
given in this case by
$$(-1)^{n(n-1)/2}=(-1)^r.
\leqno(3.4.3)$$
So (3.4.2) and hence (3.4.1) are shown.
Thus the sign of $\phi_X^n(e_1(v^{\mu}),e_1(v^{\mu}))$ coincides
with that of $\phi_X^n(u,u')$ up to the sign $(-1)^r$.

We also get the second sign $(-1)^r$ from $\gamma'$ in the
definition of $\phi_X^n$ in (3.1.2) as is shown in (3.4.3).
We then get the third sign $(-1)^r$ from the products
$$\phi_X(u_k,u'_k)\,\phi_X(u'_k,u_k)\,\,\,(k\in[1,r]),$$
since $\phi_X$ is anti-symmetric on $V^{\sb}_{\rm odd}$
and $\phi_X(u_k,u'_k)\in\R$.

Thus we get the sign $(-1)^r$ in total (since we got it three times).
This sign coincides with that of the corresponding term on the
right-hand side, which is the sign of the coefficient of $t^{2r}$
in the polynomial $(1-t^2)^m$ where $m:=\sum_{i\in\N}\rho_{2i+1}$.
The absolute value of the coefficient is $\binom{m}{r}$, and this
also coincides with that for the left-hand side in this case.
So Theorem~2 is proved.

\ms\nin
{\bf Remark~3.5.}
The above calculation in the even degree case is closely related to
[MG2].
In the odd degree case, however, we get an anti-symmetric pairing,
and this is different from [MG1].

\ms\nin
{\bf 3.6.~Abstract Hodge index theorem.}
Let $(V^{\sb};l,\phi)$ be a graded $\R$-Hodge structure of Lefschetz
type of weight $w$ in [Sa1], Sect.\ 4 with the precise signs
(see [D4]).
This means that $V^k$ is a pure $\R$-Hodge structure of weight $w+k$
endowed with a morphism of Hodge structures
$l:V^{\sb}\to V^{\sb+2}(1)$
and a graded-symmetric pairing of vector space
$\phi:V^{\sb}\otimes V^{\sb}\to \R$ such that $\phi$ induces a
self-pairing of graded $\R$-Hodge structures with value in $\R(-w)$,
we have $\phi(lu,v)=\phi(u,lv)$ for any $u,v$,
and $(-1)^{k(k-1)/2}\phi(id\otimes l^k)$ gives a polarization of
Hodge structure on the primitive part
$V^{-k}_{\rm prim}:=\Ker\,l^{k+1}\subset V^{-k}\,\,(k\in\N)$.
These conditions imply
$$i^{q-p}\phi(l^kv,\overline{v})>0\,\,\,\,\h{for any}\,\,\,
v\in V^{p,q}_{{\rm prim},\C}\setminus\{0\},$$
where $k:=w-p-q\in\N$.
Here we use the Hodge decomposition
$$V^{-k}_{{\rm prim},\C}=\mopl_{p+q=w-k}V^{p,q}_{{\rm prim},\C}.$$
In some references, $i^{p-q}$ is used instead of $i^{q-p}$.
However, this does not cause a problem if $p+q$ is even.
Set
$$\chi_y(V^{\sb})=\sum_{p,q}(-1)^{q-w}\,h^{p,q}(V^{\sb})\,y^p
\q\h{with}\q h^{p,q}(V^{\sb})=\dim_{\C}\Gr_F^pV^{p+q-w}_{\C}.$$
Since $\chi_{-y}(V^{\sb})=\sum_{p,q}(-1)^{p+q-w}\,h^{p,q}(V^{\sb})
\,y^p$, this agrees with the previous definition.

Assume $w$ is {\it even} so that $(-1)^w=1$.
Let $\s(\phi|V^0)$ denote the signature of the restriction of the
graded-symmetric self-pairing $\phi$ to $V^0$.
Then
$$\s(\phi|V^0)=\chi_1(V^{\sb})=\sum_{p,q}(-1)^q\,h^{p,q}(V^{\sb}).
\leqno(3.6.1)$$
This follows from the same calculation as in [Hi], p.~125,
Thm.~15.8.2, using the above conditions together with the primitive
decomposition.
Since $w$ is even, we also have
$$\chi(V^{\sb})=\chi_{-1}(V^{\sb})=
\sum_{p,q}(-1)^{p+q}\,h^{p,q}(V^{\sb}).
\leqno(3.6.2)$$

\ms\nin
{\bf 3.7.~Relation between Corollary~3 and Theorem~2.}
Using the above properties, we can show that Corollary~3 implies
Theorem~2 in the case when $X$ is a projective variety, the complex
$K$ is an intersection complex underlying a polarizable Hodge module
$\Mc$ of even weight, and $\phi$ is a polarization of $\Mc$.
In this case $H^{\sb}(X,\Mc)$ is a graded Hodge structure of
Lefschetz type of weight w, see [Sa1], Th.~5.3.1.
Then Theorem~2 is shown by splitting the summation inside the
exponential in the last term of the formula in Corollary~3 in two
parts according to the parity of the index of summation.
(For a similar formulation in case X is a smooth projective variety,
see [Mo], p.~173, Cor.~2.13, which is based on the Hodge index
theorem for global projective complex $V$-manifolds in
loc.~cit., p.~171, Cor.~2.11.)

\ms\nin
{\bf 3.8.~Multiple K\"unneth formula for $A$-complexes.}
In case of bounded $A$-complexes with constructible cohomology
sheaves on topological stratified spaces as in (3.1) (or [Ve]),
the multiple K\"unneth formula holds by assuming the finiteness
of the global cohomology, where $A$ a field of characteristic 0.
Indeed, let $K_i\in D^+_c(X_i,A)$ with
$\dim H^{\sb}(X_i,K_i)<\infty$ for $i\in[1,n]$.
There is a {\it canonical morphism} of complexes
$$\motim_{i=1}^n\R\Gamma(X_i,K_i)\to
\R\Gamma\bigl(\mprod_iX_i,\bt_{i=1}^nK_i\bigr),
\leqno(3.8.1)$$
which is defined by taking a flasque resolution
$$\bt_{i=1}^nK_i\simto\Kct.
\leqno(3.8.2)$$
Here we may assume that each $K_i$ is flasque by replacing it with
a flasque resolution if necessary.
It is shown that the canonical morphism (3.8.1) is a
quasi-isomorphism as follows.

Let $pr_n$ denote the projection to the $n$-th factor $X_n$ with
fiber $\X':=\mprod_{i=1}^{n-1}X_i$ where $n\ge 2$.
Note that (3.8.1) holds for $\Kc':=\bt_{i=1}^{n-1}K_i$ by inductive
assumption if $n-1>1$ (and it is trivial if $n-1=1$).
We first show the quasi-isomorphism
$$\R\Gamma(\X',\Kc')\sotim K_n\simto\R(pr_n)_*(\bt_{i=1}^nK_i),
\leqno(3.8.3)$$
where $\R\Gamma(\X',\Kc')$ on the left-hand side is identified with
a constant sheaf complex on $X_n$.
The morphism in (3.8.3) is defined by using the flasque resolution
(3.8.2) (together with the inductive hypothesis for $\Kc'$).
To show that (3.8.3) is a quasi-isomorphism, we have to determine
the stalk of the right-hand side at each $x_n\in X_n$.
For this we take $U_{x_n}\subset X_n$ in (3.1) and prove the
following canonical quasi-isomorphism induced by the restriction
morphism for the inclusion
$\X'\times\{x_n\}\into\X'\times U_{x_n}$:
$$\R\Gamma(\X'\times U_{x_n},\Kc'\bt(K_n|U_{x_n}))\simto\R\Gamma
(\X',\Kc'\sotim K_{n,x_n}).
\leqno(3.8.4)$$
Note that the right-hand side is canonically isomorphic to
$\R\Gamma(\X',\Kc')\sotim K_{n,x_n}$ using the hypotheses on the
finiteness, and the left-hand side is independent of the size of
$U_{x_n}$ under the restriction morphisms using the cone structure
in (3.1) since $K_n$ is cohomologically constructible with respect
to the stratification in (3.1).
For the proof of (3.8.4), consider the direct image by the
projection $pr'$ to $\X'$ with fiber $U_{x_n}$.
We get a canonical quasi-isomorphism induced by the restriction
under the inclusion $\X'\times\{x_n\}\into\X'\times U_{x_n}$:
$$\R pr'_*(\Kc'\bt(K_n|U_{x_n}))\simto\Kc'\sotim K_{n,x_n}.
\leqno(3.8.5)$$
This is proved by restricting it over $\prod_{i=1}^{n-1}U_{x_i}$
and reducing to the fact that the resolution (3.8.2) induces a
quasi-isomorphism at each stalk.
Then (3.8.4) and hence (3.8.3) follow by applying the global
section functor over $\X'$ to (3.8.5).
We now apply the global section functor over $X_n$ to
(3.8.3), and conclude that (3.8.1) is a quasi-isomorphism by
increasing induction on $n$.
(The argument seems to work also for $K_i\in D^b_c(X_i,A)$ in the
sense of [Ve] assuming $\dim H^{\sb}(X_i,K_i)<\infty$.)

The above quasi-isomorphism (3.8.1) implies the multiple K\"unneth
isomorphism and also a remark after Theorem~1.
Note that similar assertions hold for cohomology with compact
supports where $U_{x_i}$ is replaced by its closure in $X_i$
after shrinking it.
(See [Sc1] for an argument using the base change theorem.)

\end{document}